\documentclass[reqno]{amsart}
\usepackage{algorithm}
\usepackage{booktabs}
\usepackage{subfig}
\usepackage[english]{babel}
\usepackage{graphicx}
\usepackage[hidelinks]{hyperref}
\usepackage{amssymb}
\usepackage{amsmath}
\usepackage[foot]{amsaddr}
\usepackage{amsfonts}
\usepackage{mathtools}
\usepackage{algpseudocode}
\usepackage{algorithm}
\usepackage[top=1in, bottom=1.25in, left=1.25in, right=1.25in]{geometry}
\usepackage [autostyle, english = american]{csquotes}
\MakeOuterQuote{"}

\usepackage{verbatim}
\usepackage{mathrsfs}
\usepackage[svgnames]{xcolor}
\usepackage{dsfont}

\usepackage{MnSymbol}

\newcommand{\underbrset}[2]{\underset{#1}{\underbrace{#2}}}

\usepackage{esvect}
\usepackage{enumitem}

\usepackage{gensymb} 
\usepackage{physics}

\newtheoremstyle{remarkStyle}{\parskip}{}{}{}{\bfseries}{.}{.5em}{}
\theoremstyle{plain} \newtheorem{definition}{Definition}[section]
\theoremstyle{plain} \newtheorem{theorem}{Theorem}[section]
\theoremstyle{remarkStyle} 
\theoremstyle{remarkStyle} 
\theoremstyle{plain} 
\theoremstyle{plain} 
\theoremstyle{plain} 
\theoremstyle{remarkStyle}

\graphicspath{{images/}}
\usepackage{float}

\begin{document}
\title[MOR for bifurcating phenomena in FSI problems]{Model order reduction for bifurcating phenomena in Fluid-Structure Interaction problems}

\author{Moaad Khamlich$^{1}$}
\author{Federico Pichi$^{1,2}$}
\author{Gianluigi Rozza$^{1}$}
\address{$^1$ mathLab, Mathematics Area, SISSA, via Bonomea 265, I-34136 Trieste, Italy}
\address{$^2$ Chair of Computational Mathematics and Simulation Science, \'Ecole Polytechnique F\'ed\'erale de Lausanne, 1015 Lausanne, Switzerland}
\begin{abstract}
This work explores the development and the analysis of an efficient reduced order model for the study of a bifurcating phenomenon, known as the Coand\u{a} effect, in a multi-physics setting involving fluid and solid media.
Taking into consideration a Fluid-Structure Interaction problem, we aim at generalizing previous works towards a more reliable description of the physics involved. In particular, we provide several insights on how the introduction of an elastic structure influences the bifurcating behaviour.
We have addressed the computational burden by developing a reduced order branch-wise algorithm based on a monolithic Proper Orthogonal Decomposition. We compared different constitutive relations for the solid, and we observed that a nonlinear hyper-elastic law delays the bifurcation w.r.t.\ the standard model, while the same effect is even magnified when considering linear elastic solid.
\\

\textbf{Keywords}: reduced order modelling; parametrized fluid-structure interaction problem; monolithic method; proper orthogonal decomposition; bifurcation theory; continuum mechanics; fluid dynamics; Coand\u{a} effect.
\end{abstract}

\maketitle
\section{Introduction and motivation}
\label{sec:intro}
In the present work, we have focused on a symmetry breaking bifurcation for the solutions of a fluid-structure interaction (FSI) problem in a planar contraction-expansion channel. We have chosen to tackle this problem because it leads to complex and rich dynamics for relatively small Reynolds numbers. In fact, the model is characterized by the so-called \textit{Coand\u{a}} effect \cite{coanda-ahmed}, which in this setting causes the onset of an asymmetric flow that attaches to the nearby wall, despite the symmetry of both the geometry and the boundary conditions.

\begin{figure}[h]
    \centering
    \begin{minipage}{0.45\linewidth}
    \centering
    \includegraphics[width=0.65\linewidth]{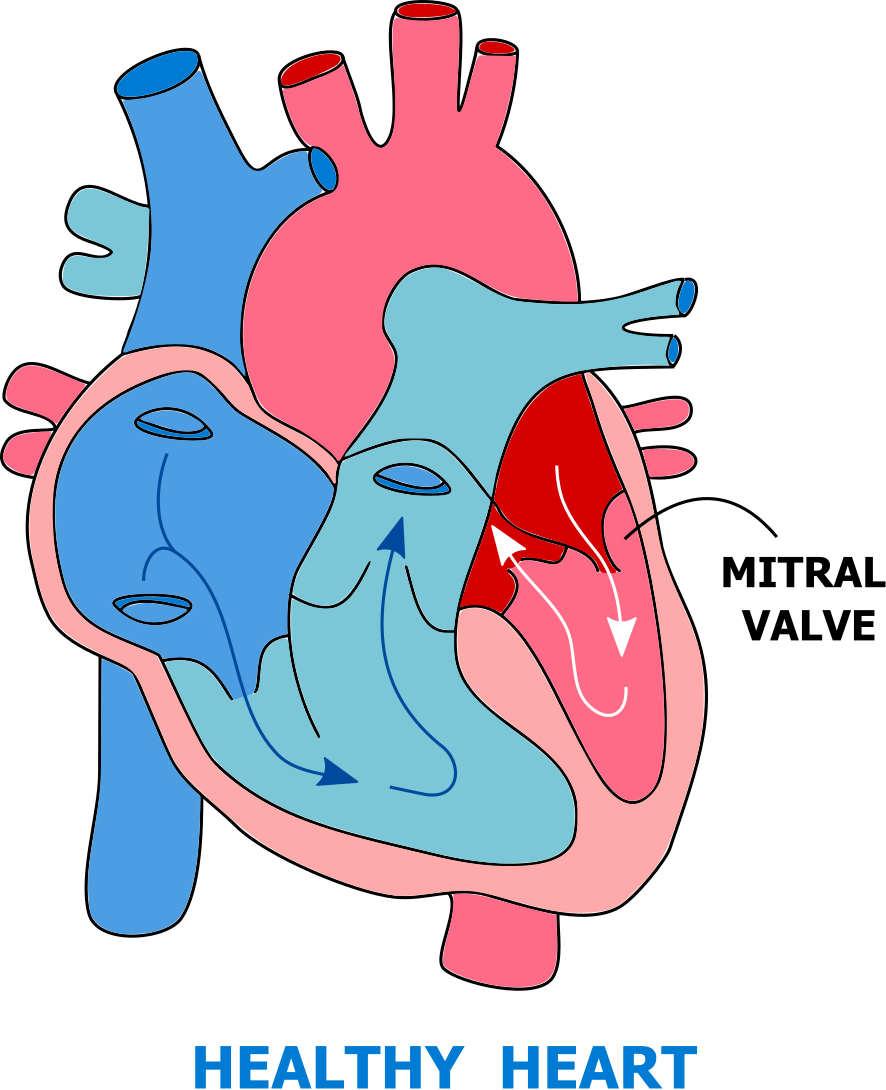}
    \end{minipage}\quad\quad
    \begin{minipage}{0.45\linewidth}
    \centering
    \includegraphics[width=0.65\linewidth]{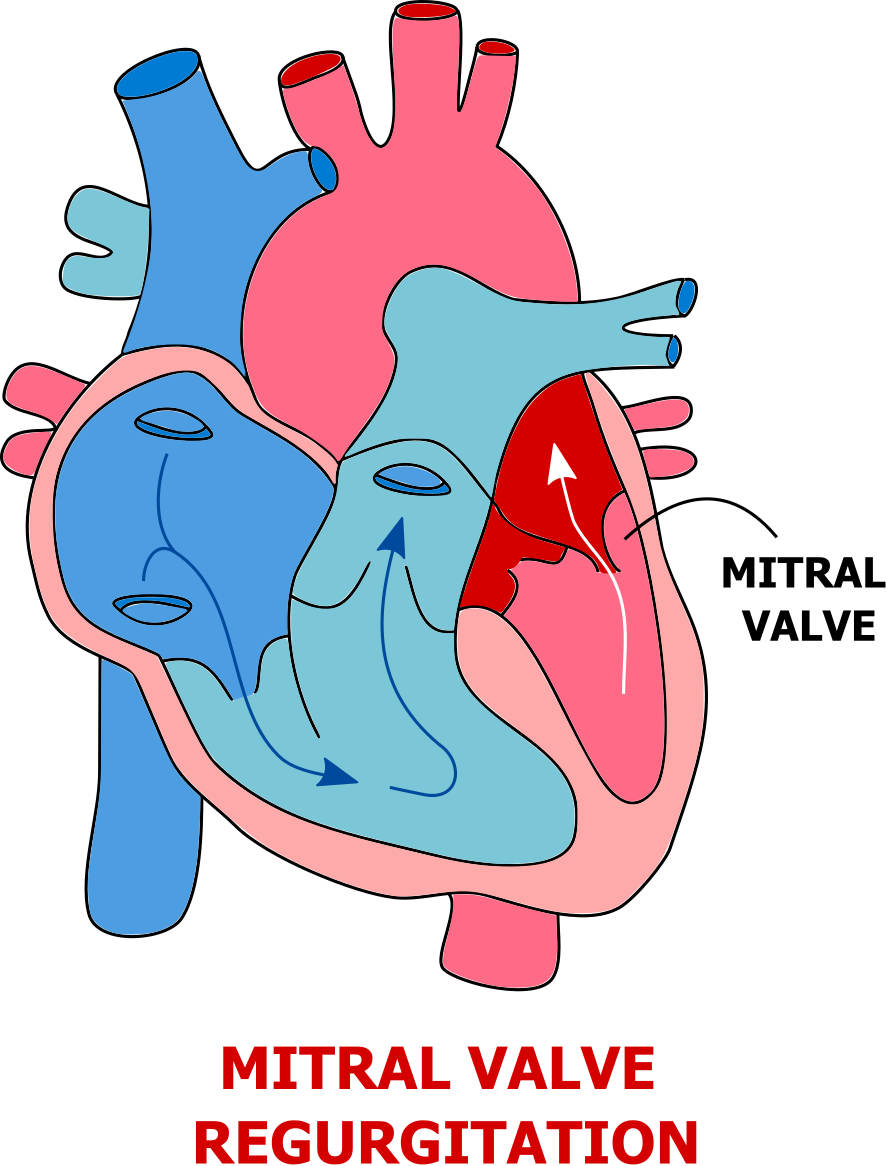}
    \end{minipage}
    \vspace*{3mm}
    \caption{Cardiac disease related to the Coanda effect: healthy recirculation (left) and mitral valve regurgitation (right).}%
    \label{fig:regurgitation}
\end{figure}

In particular, the application that motivated our analysis is a heart disease known as mitral valve regurgitation, which occurs when blood flows abnormally from the left ventricle to the left atrium because of the defective closure of the mitral valve (see Figure \ref{fig:regurgitation}). In particular situations, this anomalous flow also exhibits the aforementioned wall-hugging behavior, leading to incorrect readings of the blood volume through echocardiography, which is used to assess the severity of this pathology \cite{echoassestment}. Since a channel with a sudden expansion preserves the essential properties of the blood flow through the two heart chambers, it is well suited to characterize the main features of this phenomenon.
In fact, the test case under consideration undergoes a bifurcation due to the non-linearity of the Navier-Stokes system, and this results in a change of both the number of equilibrium solutions and their nature for a small variation of the Reynolds number $Re$.
For low $Re$,  the symmetric solution is unique and stable \cite{Serrin,Shapira, Sobey}; however, as $Re$ increases, the inertial effects become progressively more important, breaking the symmetry of the solution for a certain critical value $Re^*$ \cite{Battaglia}. Below this critical value, a physically unstable configuration of symmetric flow and two stable asymmetric ones coexist. The existence of multiple solutions for the same parameter values, which is typical of branching phenomena, gives rise to a bifurcation pattern known as \textit{supercritical pitchfork bifurcation}.

The Coand\u{a} effect has been widely studied in the context of pure fluid dynamics (see Section \ref{sec:related}); thus, we propose a generalization of the previous works by facing this phenomenon in a multi-physics context. To provide a simplified description of mitral regurgitation, we set up the FSI problem, which deals with the mutual inﬂuence between a deformable structure and the fluid flow around it. The effects of this interaction can be significant or even predominant due to the stress field exchanged at the interface.
To this end, we used the Arbitrary Lagrangian-Eulerian formulation \cite{aleformulation}, which is needed to solve the problem on a fixed reference configuration.
Our investigation aims at understanding how the introduction of the structure can modify the bifurcating behavior. For this purpose, we carried out a comparison between the multi-physics case and the one where a non-deformable rigid body replaces the elastic structure. We further deepened our investigation by considering large scale deformations.
Moreover, towards a better description of the physics involved, we introduced a nonlinear constitutive relation for the structure, i.e. the Saint Venant-Kirchoff model.

The study of the bifurcation implies being able to obtain solutions that share the same qualitative properties (which are said to belong to the same \textit{branch}) while changing the value of the parameter, which in our case will be the kinematic viscosity of the fluid, strictly linked to the aforementioned Reynolds number. This task is not trivial since, after the bifurcation has occurred, we lose the uniqueness of the solution \cite{seydel}. Consequently, we employed an algorithm combining the Galerkin-Finite Element method \cite{quarteroni2017numerical}, the Newton method \cite{ciarlet2013linear}, and a continuation method \cite{continuation}.
This approach allows us to obtain a sequence of solutions associated with the same prescribed branch.

Despite this, since the high-fidelity models usually require unaffordable computational costs due to fine mesh discretization and repeated evaluations, we also focused on the application of Reduced Order Models (ROM) \cite{benner2017model, morepas2017}. In particular, we employed the \textit{Reduced Basis} (RB) method \cite{hesthaven2015certified, patera07:book}, which has been already used to trace the bifurcating behavior of many problems in fluid-dynamics \cite{HESS2019379,pichistrazzullo,Pitton_2017}, with further applications to other fields such as quantum mechanics \cite{pichiquaini_old}, computational mechanics \cite{pichirozza} and bio-mathematics \cite{coanda}.
This method is based on projecting the nonlinear governing equations on a reduced space of much lower dimension. The latter is generated by a small set of RB functions, which can be obtained employing the Proper Orthogonal Decomposition (POD) \cite{POD}.

Our results suggest that although the bifurcation is caused by the nonlinearity of the constitutive equation for the fluid, the presence of the solid and the particular constitutive relation used to model it, scale the bifurcating behavior.
Moreover, when dealing with large scale deformations, the bifurcation point translates into the parameter interval. This effect was observed for both the linear and nonlinear solid constitutive relation; however, it is amplified in the first case.

The work is organized as follows:
\begin{itemize}
    \item   Section \ref{sec:related} provides the state of the art related to the Coand\u{a} effect in planar contraction-expansion channels;
    \item Section \ref{sec:problem_formulation} introduces the FSI problem, providing the variational formulation for the steady case;
    \item Section \ref{sec:numerical_approximation} describes all the ingredients to tackle numerically a generic bifurcation problem held by a nonlinear PDE: the high order approximation, the Reduced Basis method, and the Newton method.
    \item Section \ref{sec:numerical} shows the results for the FSI problem concerning the rigid body model, the linear elastic model, and the Saint Venant-Kirchoff nonlinear model.
\end{itemize}

\section{Related Work}\label{sec:related}
As anticipated in the introduction, the object of our investigation is a symmetry-breaking bifurcation for a planar incompressible flow through a contraction-expansion channel.

The appearance of an asymmetric behavior for laminar flows, undergoing a sudden section change, is a problem that has been widely studied both experimentally and numerically \cite{Cherdron, Sobey, Tritton}. Such phenomenon is caused by the Coand\u{a} effect, which concerns velocity fluctuations that lead to a transversal pressure gradient capable of maintaining the flow asymmetry \cite{Tritton, Wille}.

Cherdron et al. \cite{Cherdron} presented a first large experimental analysis that uses flow visualization and laser-Doppler anemometry. This study explains the origin of the phenomenon through the disturbances generated at the corners of the expansion, and subsequently amplified by shear layers. The latter plays a fundamental role, as the Moffat eddies \cite{Moffatt} forming inside them are interdependent due to confinement in the transversal direction.

Another important contribution was provided by Sobey \& Drazin \cite{Sobey}, who dealt with the problem using an asymptotic, numerical and experimental approach. They understood the strong connection of the problem with bifurcation theory, being able to deduce the presence of a Hopf bifurcation characterized by the existence of a time-periodic solution, whose $Re^{*}$ have been successively estimated by the work of Quaini et al. \cite{Quaini}. The very nature of the temporal dependence of the solution was then investigated by Fearn et al. \cite{Fearn}, who managed to trace the phenomenon to three-dimensional effects through an experimental investigation.

The results of all these studies highlighted the onset of the Coand\u{a} Effect for relatively low Reynolds numbers. This can be partially justified by the fact that although the Coand\u{a} effect occurs more strongly in turbulent flows, it is amplified in the two-dimensional setting \cite{Wille}. However, these early works were carried out working with fixed expansion ratio $\lambda$, or with few values. Battaglia et al. \cite{Battaglia} and Drikakis \cite{Drikakis} then generalized the previous discussion by studying the relation between $Re^{*}$ and $\lambda$. Their results showed that the critical Reynold number increases as the expansion ratio decreases; however, subsequent investigations, that include a more complex case study, show a non-trivial behavior of the bifurcation with respect to the expansion ratio \cite{pichi2021artificial}

The physical interpretation of the laminar flow's dynamics inside the contraction-expansion channel was provided by Hawa \& Rusak's work \cite{Hawa}, who have shown that, for low Reynolds numbers, the flow is stabilized by viscous dissipation. As $Re$ increases, this stabilization is not sufficient to counteract the destabilization by the convective effects upstream of the expansion, and precisely this imbalance leads to the instability of the symmetric solution.

Subsequent investigations mainly focused on different kinds of generalizations, such as the constitutive relationship \cite{Mishra}, geometry \cite{Mizushima}, or the link between numerical schemes and the critical Reynolds value \cite{Quaini}.

The analogy between the flow in the contraction-expansion channel and the anomalous jet due to mitral regurgitation (see Section \ref{sec:intro}), have prompted recent works \cite{Pichi,coanda,Pitton_2017} to create a connection between the medical literature \cite{Chandra,Little,Vermeulen} on this pathology and the numerical modeling of the cardiovascular system. In particular, these works have tackled the problem of obtaining efficient and real-time simulations, using reduced order models, to be exploited in a real-life scenario. For example, the work of Pitton et al.  \cite{coanda} allowed to reproduce the Coand\u{a} effect in a mock heart chamber.

The big limitation of these previous works lies in the fact that they have faced the problem by solving the Navier-Stokes system in the contraction-expansion channel with fixed straight walls. An important step towards more realistic test cases was made by Hess et al. \cite{quaini-hess}, by studying a contraction-expansion channel with curved walls.

Our work seeks to carry out this attempt at generalization. In fact we are not aware of studies on the effect that the coupling of the fluid with an elastic solid can have on the bifurcation. Therefore we have decided to tackle the problem in the context of fluid-structure interaction.
In addition to providing a replicable benchmark for future investigation in the multi-physics context, we tested a reconstruction using the RB method, which has been proven to be effective for FSI problems \cite{Bal2016, Ballarin2017,Nonino, fluids6060229}.

\section{Problem formulation}%
\label{sec:problem_formulation}
Let us consider the abstract strong form of a parametric PDE, which reads: given a parameter $ \boldsymbol {\mu} \in \mathbb {P}, $ find $ X (\boldsymbol {\mu}) \in \mathbb{X}$ such that
\begin {equation}
    \label {strongfm}
    F (X (\boldsymbol {\mu}); \boldsymbol {\mu}) = 0, \quad F: \mathbb {X} \times \mathbb {P} \to \mathbb {X} ^ {\prime},
\end {equation}
being $ \mathbb {X} $ the functional space where we are seeking the solution $ X $, and $ \mathbb {X}' $ its dual space.
In order to retrieve a variational formulation of the problem, we introduce the following variational form:
$$
f (X(\boldsymbol{\mu}), Y;\boldsymbol{\mu}) \doteq \langle F (X(\boldsymbol{\mu}); \boldsymbol{\mu}), Y \rangle, \quad f: \mathbb {X} \times \mathbb {X} \rightarrow \mathbb {R},
$$
where $Y \in \mathbb{X}$ is a test function and $\langle \cdot, \cdot \rangle$ denotes the duality pairing between $\mathbb{X}$ and $\mathbb{X}^{ \prime}$.

Then, the weak formulation of equation \eqref {strongfm} reads: given $ \boldsymbol {\mu} \in \mathbb {P}, $ find $ X (\boldsymbol {\mu}) \in \mathbb{X}$ such that
\begin {equation}
    \label {variationalform}
    f (X(\boldsymbol{\mu}), Y; \boldsymbol{\mu}) = 0 \quad \forall\ Y \in \mathbb {X}.
\end {equation}

We are interested in studying the behavior of the solution as the value of the parameter varies. In particular, if we assume that $F$ is continuously differentiable, and we denote by $D_{X} F(X; \boldsymbol{\mu}): \mathbb{X}\to \mathbb{X}^{\prime}$ its Fr\'echet derivative with respect to $X$, then the following important result holds true \cite{ma2005bifurcation}:

\begin{theorem}
    Let $(\bar{X},\bar{\boldsymbol{\mu}})$ be a known solution of equation \eqref{strongfm}.  Let $B_{r}(\bar{X}), B_{\bar{r}}(\bar{\boldsymbol{\mu}})$ be two balls of radius $r$ and $\bar{r}$ around $\bar{X}$ and $\bar{\boldsymbol{\mu}}$, respectively. Let $F: \mathbb{X} \times \mathbb{P} \rightarrow \mathbb{X}^{\prime}$ be a $C^{1}$ map.
    If  $D_{X} F(\bar{X} ; \bar{\boldsymbol{\mu}})$ is bijective
    then, there exist $r, \bar{r}>0$ and a unique solution $X(\boldsymbol{\mu}) \in B_{r}(\bar{X}) \cap \mathbb{X}$ such that
$$
F(X(\boldsymbol{\mu}) ; \boldsymbol{\mu})=0 \quad \forall \boldsymbol{\mu} \in B_{\bar{r}}(\overline{\boldsymbol{\mu}}) \cap \mathbb{P}.
$$
\end{theorem}
Therefore, when the aforementioned assumptions fail, it can happen that even a slight a variation of the parameter leads to changes in the number of the existing solutions for a given value of the parameter \cite{caloz1997numerical,Brezzi1980}.
In particular, this means that multiple solutions could branch off from a reference solution $ \tilde{X} (\boldsymbol{\mu}) $.
The value of the parameters for which this happens will be called \textit{bifurcation points} of equation \eqref{strongfm}, and can be defined as done in \cite{ambrosetti_malchiodi_2007}:
\begin{definition}
    \label{defbif1}
    A parameter value $\boldsymbol{\mu}^{*} \in \mathbb{P}$ is a bifurcation point for equation \eqref{strongfm} if there exists a sequence $\left(X_{n}, \boldsymbol{\mu}_{n}\right) \in \mathbb{X} \times \mathbb{P},$ with $X_{n} \neq \tilde{X}(\boldsymbol{ \mu _{n}}) $ such that:
\begin{itemize}
    \item $F\left(X_{n} ; \boldsymbol{\mu}_{n}\right)=0$;
    \item $\left(X_{n}, \boldsymbol{\mu}_{n}\right) \rightarrow\left(\tilde{X}(\boldsymbol{\mu^{*})} , \boldsymbol{\mu}^{*}\right)$.
\end{itemize}
\end{definition}

The manifold represented by the solutions with the same qualitative behavior will be named \textit{branch} and denoted by $ \mathcal{M} $. The objective of our investigation is the reconstruction of the entire solution manifold, given by the union of all the branches:
\begin{equation}
    \label{solutionmanif}
\mathcal{S}=\bigcup_{i=1}^{k} \mathcal{M}_{i}:=\bigcup_{i=1}^{k}\left\{X_{i}(\boldsymbol{\mu}) \in \mathbb{X} \mid F\left(X_{i}(\boldsymbol{\mu}) ; \boldsymbol{\mu}\right)=0, \boldsymbol{\mu} \in \mathbb{P}\right\}.
\end{equation}
In particular, this is fundamental when trying to recover the so-called \textit{bifurcation diagram}. The latter is an important tool, typical of the study of dynamical systems, which allows us to investigate the qualitative behavior of the solutions as a function of the parameters' set taken into consideration. The idea is to trace the evolution of a significant scalar quantity $ s (X (\mu)) $, which can be considered as the output of the dynamical system. Typical examples in this regard are the norm of the solutions, or a pointwise evaluation of one of the solution's components. The presence of multiple solutions corresponding to the same parameter value results in multiple branches in the bifurcation diagram.

In the following, we want to apply the former discussion to the fluid-structure interaction problem. This will be done by retrieving both the strong and weak form for the particular problem at hand.
However, we need to first introduce the \textit{Arbitrary Lagrangian-Eulerian} formulation, which will be used to solve the problem encoding coherently the fluid's flow and the displacement of the structure.

\subsection{Preliminaries: configurations and ALE description}
\label{sec:Ale_formulation}

%

In this section, we are going to introduce the Arbitrary Lagrangian-Eulerian (ALE) formulation \cite{aleformulation} which is widely used for the description of the fluid motion in FSI problems \cite{fsi-book}. In particular, let $\Omega(t) $  be the physical domain, which in our case will be a subset of the two-dimensional Euclidean space $\mathbb{R}^{2}$.
We assume that the solid resolution domain $\Omega_{s}(t)$ and the fluid one $\Omega_{f}(t)$ are such that  $\Omega(t)=\Omega_{f}(t) \cup \Omega_{s}(t)$ and they do not overlap\footnote{In the following, we will use the subscript $_{s}$ (respectively $_{f}$) to refer to \textit{solid} (respectively \textit{fluid}) 's phase quantities.}.

\begin{figure} [htbp]
\centering
\def\svgwidth{0.85\linewidth}
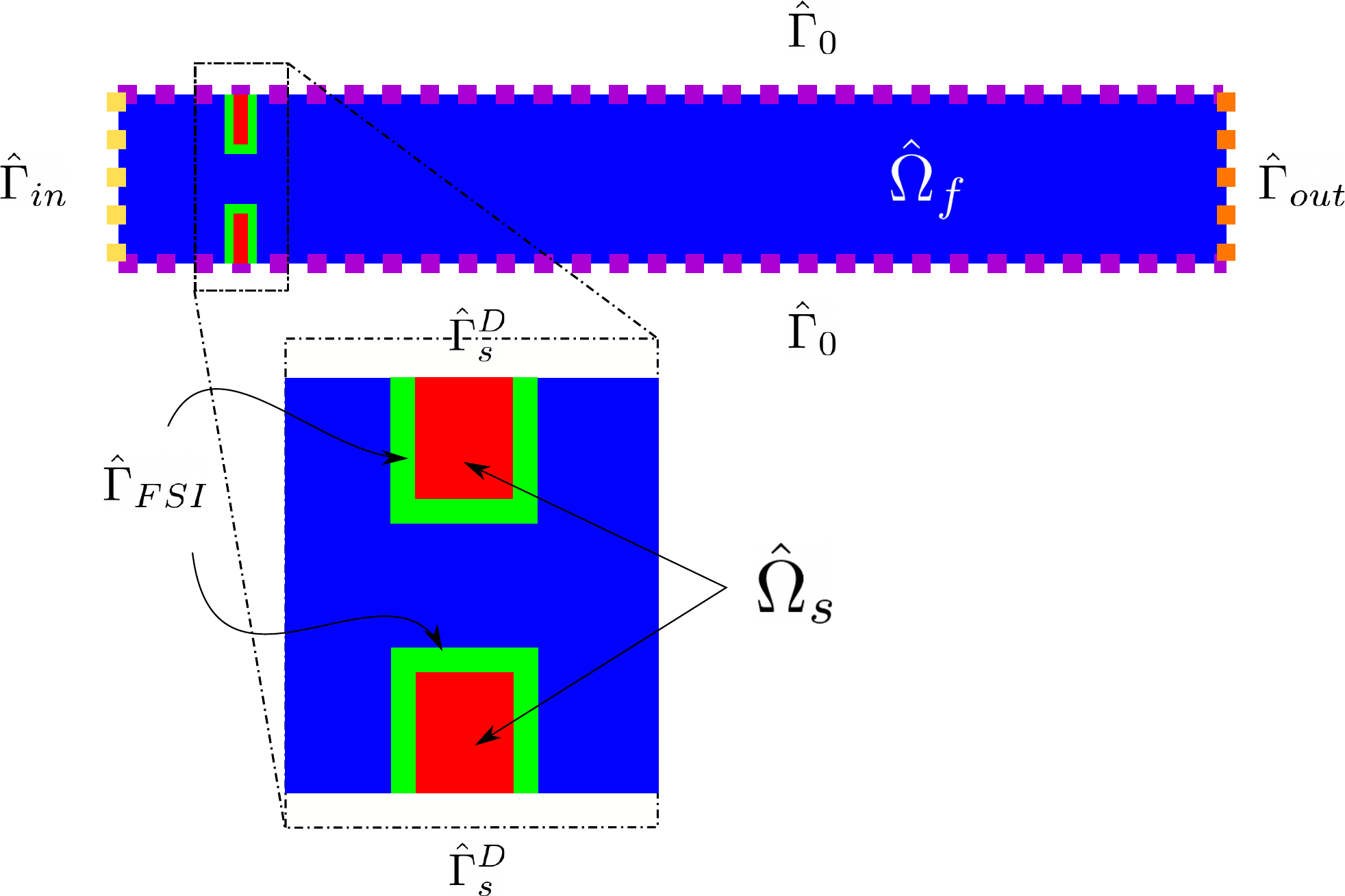
\caption{Reference configuration for the FSI problem.}
\label{fig:domain}
\end{figure}

We know that for a generic FSI problem we have to deal with a moving domain. As regards the description of the solid phase, the problem is solved by defining the equations in material form, i.e.\ in the reference domain $ \hat {\Omega} _ {s} $. Typically the problem does not concern fluids as they adapt to the domain by assuming its form. However, in fluid-structure interaction problems, the fluid resolution mesh needs to be in motion in order to be consistent with the structural deformation. One possible approach to deal with this issue is to use the already mentioned ALE formulation, which consists in defining an intermediate reference system $ \hat{\Omega}_{f}$ on which the fluid equations are pulled back.
Therefore, it is necessary to introduce an additional field $ \boldsymbol{d}_f $, which represents the displacement of the fluid mesh we are modeling; this displacement must necessarily satisfy a Dirichlet condition on the boundary with the solid resolution domain $ \Gamma_{FSI} = \Omega_{f} (t) \cap \Omega_{s} (t) $:

$$
\boldsymbol{d}_f = \boldsymbol{d}_{s} \text{ on }  \Gamma_{FSI} \subset \partial \Omega_{f}(t);
$$
being $\boldsymbol{d}_{s}$ the solid displacement.
The crucial task is to extend the displacement $\eval{\boldsymbol{d}_{f}}_{\Gamma_{FSI}}$ inside the domain $\Omega_{f}(t)$.
To do so, as done in \cite{aleformulation}, we introduce a fictitious displacement which follows the Dirichlet datum prescribed on the boundary. This extension is said to be arbitrary because it only concerns the reference system used to describe the fluid motion, and not the motion itself. This led us to consider an extension problem, among which we have decided to use the \textit{harmonic extension}. The latter consists in solving a Poisson problem with inhomogeneous Dirichlet boundary conditions on the fluid-structure interface:
\begin{equation}\label{eq:domainHarmonicExtension}\begin{dcases}
    -\Delta \hat{\boldsymbol{d}}_f= 0 & \text{in $\hat{\Omega}_f$}, \\
    \hat{\boldsymbol{d}}_f = \hat{\boldsymbol{d}}_{s} & \text{on $\hat\Gamma_{FSI}$},\\
    \hat{\boldsymbol{d}}_f = 0 & \text{on $\hat\Gamma = \partial\hat\Omega_f\backslash\hat\Gamma_{FSI}$}.
\end{dcases}
\end{equation}

The harmonic extension has been formulated in the fixed reference configuration $\hat{\Omega}_f$, which is the inverse image of the Eulerian fluid resolution domain $\Omega_{f}(t)$ through the ALE map $\mathcal{A}_{f}(t): \hat{\Omega}_{f} \to \Omega_{f}(t)$:
\begin{equation}
\label{eq:ALEdomainUpdate}
\Omega_f(t) = \{\boldsymbol{x}= \mathcal{A}_{f}( \hat{\boldsymbol{x}},t)=\hat{\boldsymbol{x}} + \hat{\boldsymbol{d}}_f(\hat{\boldsymbol{x}},t) , \; \hat{\boldsymbol{x}}\in\hat\Omega_{f} \}.
\end{equation}

The time derivative in such a system, denoted with $\fdv{t}$, can be computed for any scalar or vector quantity $v$ through the Reynolds transport formula:
\begin{equation}
    \label{eq:eqIndEq}
    \underbrset{ \mathclap{ \text{ALE time derivative} } }{ \fdv{v}{t}} = \pdv{v}{t} + \mathbf{u_{f}}\cdot\grad v, \qquad \mathbf{u_{f}} = \frac{\partial \mathcal{A}_{f}}{\partial t}\circ \mathcal{A} _{f}^{-1}.
\end{equation}

\subsection{Strong formulation for the FSI problem}%
\label{strongformulationFSI}
In this section, we are going to introduce the fully coupled formulation for the FSI problem.
Since in this setting the ALE description is used for the fluid and the  Lagrangian one for the solid, we will need two distinct reference domains, namely one for the solid $ \hat{\Omega}_{s}$ and one for the fluid $ \hat{\Omega}_{f}$ (see Figure \ref{fig:domain} for a visual representation).
The unknowns in this setting are: the global displacement $\boldsymbol{d} =(\boldsymbol{d} _{s},\boldsymbol{d} _{f})$, the fluid velocity $ \mathbf{u}_{f}$, and the pressure $p$.

By combining the equations for the solid phase (momentum equation), the fluid phase (Navier-Stokes system) and the mesh displacement (harmonic extension) in a single system, we obtain the full formulation for the FSI problem, namely:
\begin{equation}
    \label{stfm_1}
    \begin{dcases}
        \rho_{f}\left( \fdv{\mathbf{u}_{f}}{t}  +\left( \mathbf{u}_{f}- \fdv{\boldsymbol{d_{f}}}{t}   \right) \cdot \nabla \mathbf{u}_{f} \right) -\operatorname{div}\boldsymbol{\sigma_{f}} (\mathbf{u}_{f}, p_{f}) =\boldsymbol{b} _{f} \quad &\text{in} \quad \Omega_{f}(t) \times (0,T],\\
        \operatorname{div} \mathbf{u}_{f}=0 \quad & \text{in} \quad\Omega _{f}(t) \times (0,T], \\
        -\hat{\Delta} \boldsymbol{\hat{d}}_{f}=0  \quad &\text{in} \quad \hat{\Omega}_{f} \times (0,T],\\
        \rho _{s} \frac{\partial ^{2} \hat{\boldsymbol{d} _{s}}}{\partial t ^{2} } -\hat{\operatorname{div}} \mathbf{P} (\hat{\boldsymbol{d} _{s}}) =\boldsymbol{b} _{s} \quad &\text{in} \quad \hat{\Omega}_{s} \times (0,T],
    \end{dcases}
\end{equation}
where $\rho_{s}$ and $\rho_{f}$ are respectively the solid and fluid densities, whereas $\boldsymbol{b}_{s}$ and $\boldsymbol{b}_{f}$ represent the external volume forces for the two physics.
Moreover, we denoted with $\boldsymbol{\sigma_{f}}$  the \textit{Cauchy stress tensor} and with $\mathbf{P}$ the \textit{second Piola tensor}; their expressions are provided by means of the constitutive relations for Newtonian incompressible fluid and linear elastic solids, which read:
        \begin{equation}
            \begin{dcases}
            \boldsymbol{\sigma} _{f}(\mathbf{u}_{f},p_{f})=\rho_{f} \mu \left( \nabla \mathbf{u}_{f}+ \nabla^{T} \mathbf{u}_{f} \right) -p \mathbf{I}, \\
            \mathbf{P}(\hat{\boldsymbol{d} _{s}})= \lambda_{s}\operatorname{tr} \mathbf{E}_{s}(\boldsymbol{\hat{d}} _{s}) \mathbf{I} + 2 \mu_{s} \mathbf{E}_{s}(\boldsymbol{\hat{d}} _{s}).
            \end{dcases}
        \end{equation}
$\lambda_{s}$ and $\mu_{s}$ are a couple of material-dependent quantities for the solid known as \textit{Lamè parameters}, whereas $\mathbf{E}_{s}$ is the linearized strain operator, which is defined as
\begin{equation}
    \mathbf{E}_{s}(\hat{ \boldsymbol{d}}_{s}) \doteq \frac{1}{2} ( \hat{ \nabla} \hat{\boldsymbol{d}}_{s}+\hat{ \nabla}^{T } \hat{\boldsymbol{d}}_{s}).
\end{equation}

$\mu$ is instead a fluid-dependent quantity, known as kinematic viscosity, which assumes a fundamental role in our analysis since it will be the \textit{bifurcation parameter} under investigation.

In Equation \eqref{stfm_1} we have introduced a notation clarifying the coordinate system we use in each equation: the divergence in the momentum equation for the solid reads $ \hat{\operatorname{div}}$, because it needs to be taken w.r.t.\ the reference variable $ \hat{\boldsymbol{x}}$, whereas the time derivative in the Navier-Stokes equation is referred to as $ \fdv{}{t}$ because it is the ALE time derivative which has been already introduced in Equation \eqref{eq:eqIndEq}. The same notation is used also for fields, therefore $\boldsymbol{d}_{f}$ and $\hat{\boldsymbol{d}}_{f}$ describe the same field, but the difference is in their dependence as a function of $(\hat{\boldsymbol{x}},t)$ or $(\boldsymbol{x},t)$ \footnote{It should be clear that we are considering the same field but expressed in different variables. Exactly because of this, one can obtain a \textit{ALE} description of an \textit{Eulerian} field or vice versa, using composition with $ \mathcal{A}_{f} $ or with $ \mathcal{A}_{f}^{-1} $.
\begin{equation*}
    \Psi(\boldsymbol{x} , t):=\hat{\Psi}( \mathcal{A}_{f} (\boldsymbol{x} , t), t), \quad \hat{\Psi}(\hat{\boldsymbol{x} }, t):=\Psi( \mathcal{A}_{f}^{-1}(\hat{\boldsymbol{x} }, t), t).
\end{equation*}
}.

There is still a problem with the former set of equations. In fact, even if the Eulerian time derivative has been substituted by the ALE one in the Navier-Stokes system, we have not pulled back those equations on the reference domain $ \hat{\Omega}_{f}$.  In order to do so, we introduce:
\begin{equation}
    \mathbf{F}:= \hat{\nabla} \mathcal{A}_{f} \qquad \text{and} \qquad J:=\operatorname{det} \mathbf{F};
\end{equation}
which are the gradient of the ALE map (see Equation \eqref{eq:ALEdomainUpdate}) and its determinant. We can perform a change of variables by exploiting the Piola transformation rule to express the fluid equations on the reference configuration, which leads to:
\begin{equation}
\label{fullformref}
    \begin{dcases}
        \rho_{f}J\left( \fdv{\hat{\mathbf{u}}_{f}}{t}  + \hat{ \nabla}\hat{\mathbf{u}}_{f}\mathbf{F}^{-1}\left( \hat{\mathbf{u}}_{f}- \fdv{\boldsymbol{\hat{d}_{f}}}{t}   \right)  \right) -\hat{\operatorname{div}}(J \hat{\boldsymbol{\sigma}}_{f} (\hat{\mathbf{u}}_{f}, \hat{p}_{f})\mathbf{F}^{-T}) =J\boldsymbol{\hat{b}} _{f} & \text{in} \quad \hat{\Omega}_{f} \times (0,T],\\
        \hat{\operatorname{div}} (J \mathbf{F}^{-1}\mathbf{\hat{u}}_{f})=0 & \text{in} \quad \hat{\Omega} _{f}\times (0,T], \\
        -\hat{\Delta} \boldsymbol{\hat{d}}_{f}=0  & \text{in} \quad \hat{\Omega}_{f} \times (0,T],\\
        \rho _{s} \frac{\partial ^{2} \hat{\boldsymbol{d} _{s}}}{\partial t ^{2} } -\hat{\operatorname{div}} \mathbf{P} (\hat{\boldsymbol{d} _{s}}) =\boldsymbol{b} _{s} & \text{in} \quad \hat{\Omega}_{s} \times (0,T],
    \end{dcases}
\end{equation}
with $\hat{\boldsymbol{\sigma}}_{f} $ being the representation of $\boldsymbol{\sigma}_{f}$ in the reference configuration, namely:
\begin{equation}
    \hat{\boldsymbol{\sigma}}_{f}=\mu(\hat{\nabla}\hat{\mathbf{u}}_{f}\mathbf{F}^{-1}+\mathbf{F}^{-T} \hat{\nabla}^{T}\hat{\mathbf{u}}_{f}) -\hat{p}\boldsymbol{I}.
\end{equation}

In order to deal with a well-posed problem, we need to specify some proper coupling and boundary conditions. As concerns the former, we supplemented system \eqref{stfm_1} by means of the following conditions on $\hat{\Gamma}_{FSI}$:
            \begin{equation}
                \begin{dcases}
                    \hat{\boldsymbol{d}}_{f}= \hat{\boldsymbol{d} }_{s} & \text{on} \quad \hat{\Gamma}_{FSI}, \\
                    \hat{\mathbf{u}}_{f}= \frac{\partial \hat{ \boldsymbol{d}} _{s}}{\partial t} & \text{on} \quad \hat{\Gamma}_{FSI}, \\
                    \mathbf {P} \hat {\boldsymbol{n} } = (J_{f} \boldsymbol {\hat{\sigma}}_{f} \mathbf {F}_{fluid} ^ {- T}) \hat {\boldsymbol{n} }  & \text{on} \quad \hat{\Gamma}_{FSI}.
                \end{dcases}
            \end{equation}
While for the latter we imposed the following:
\begin{equation}
\begin{dcases}
J \hat{\sigma}_{f}\left(\hat{\boldsymbol{u}}_{f}, \hat{p}_{f}\right) \mathbf{F}^{-T} \hat{\boldsymbol{n}}=-p_{i n} \hat{\boldsymbol{n}} & \text { on }   \quad \hat{\Gamma}_{i n}, \\
J \hat{\sigma}_{f}\left(\hat{\boldsymbol{u}}_{f}, \hat{p}_{f}\right) \mathbf{F}^{-T} \hat{\boldsymbol{n}}=-p_{\text {out }} \hat{\boldsymbol{n}} & \text { on }   \quad \hat{\Gamma}_{\text {out }}, \\
\hat{\boldsymbol{u}}_{f}=0 & \text { on } \quad \hat{\Gamma}_0, \\
\hat{\boldsymbol{d}}_{f}=0 & \text { on }   \quad\hat{\Gamma}_0,\\
\hat{\boldsymbol{d}}_{s}=0 & \text { on }   \quad\hat{\Gamma}_{s}^{D}.
\end{dcases}
\end{equation}
In the previous equations $\hat{\boldsymbol{n}}$ is the normal vector outgoing from $\hat{\Omega}_{f}$. See Figure \ref{fig:domain} for a sketch of the different portions of the boundaries.

As for now, everything is expressed in the reference configuration $ \hat{\Omega}=\hat{\Omega}_{f} \cup \hat{\Omega}_{s}$, therefore, since there is no need to distinguish between different coordinate systems, we will simplify the notation removing the hat symbol.

\subsection{Variational formulation for the steady case}%
\label{sub:variational_formulation}
We proceed now recovering the variational formulation of the problem \eqref{fullformref}, which is needed in order to solve the system by means of a Finite Elements discretization.
First of all, since the model will be solved in a steady setting, let us proceed to recover its strong form.  From the solid's momentum equation, we observe that the acceleration in the reference domain $\frac{\partial ^{2} \hat{\boldsymbol{d} _{s}}}{\partial t ^{2} } $  is equal to 0.
Whereas, for what concerns the Navier-Stokes system, one can proceed by annihilating the partial derivative in time for the Eulerian framework before doing the pull-back on the reference configuration. This leads to the following (steady) system:
\begin{equation}
    \label{steadyfsi}
    \begin{dcases}
        \rho_{f}J  \nabla\mathbf{u}_{f}\mathbf{F}^{-1}\mathbf{u}_{f}   -\operatorname{div}(J \boldsymbol{\sigma}_{f} (\mathbf{u}_{f}, p_{f})\mathbf{F}^{-T}) =J\boldsymbol{b} _{f} \quad &\text{in} \quad \Omega_{f},\\
        \operatorname{div} (J \mathbf{F}^{-1}\mathbf{u}_{f})=0 \quad &\text{in} \quad \Omega _{f}, \\
        \Delta \boldsymbol{d}_{f}=0  \quad &\text{in} \quad \Omega_{f},\\
        -\operatorname{div} \mathbf{P} (\boldsymbol{d} _{s}) =\boldsymbol{b} _{s} \quad &\text{in} \quad \Omega_{s}.
    \end{dcases}
\end{equation}
Furthermore, also the coupling condition needs to be changed accordingly:

            \begin{equation}
                \begin{dcases}
                    \boldsymbol{d}_{f}= \boldsymbol{d} _{s} \quad &\text{on} \quad \hat{\Gamma}_{FSI}, \\
                    \mathbf{u}_{f}= 0 \quad &\text{on} \quad \hat{\Gamma}_{FSI}, \\
                    \mathbf {P}   {\boldsymbol{n} } = (J \boldsymbol {\sigma}_{f} \mathbf {F}_{fluid} ^ {- T})  {\boldsymbol{n} }  \quad & \text{on} \quad \hat{\Gamma}_{FSI}.
                \end{dcases}
            \end{equation}

Finally, if we define the functional setting as $V_{f}=E_{f}=H^{1}(\hat{\Omega}_{f};\mathbb{R}^{2}), Q_{f}=L^{2}(\hat{\Omega}_{f})$ and $V_{s}=H^{1}(\hat{\Omega}_{s};\mathbb{R}^{2})$, then we can multiply each equation of system \eqref{steadyfsi} by an appropriate test function belonging to the corresponding functional space and apply the divergence theorem. This leads to the following weak formulation:
\begin{itemize}
    \item  \textbf{FLUID-STRUCTURE MOMENTUM EQUATION:}
\begin{equation}
    \label{fluidstructmomentum}
    \begin{aligned}
    &+\int _{ \hat{\Omega}_{f}} J \rho_{f} [ \nabla \mathbf{u}_{f} \mathbf{F}^{-1}]\mathbf{u}_{f}\cdot \mathbf{v}_{f}\ d \hat{x} + \int _{ \hat{\Omega}_{f}} J \boldsymbol{\sigma}_{f} \mathbf{F}^{-T}: \nabla \mathbf{v}_{f} \ d \hat{x} \\
    &+\int_{\hat{\Gamma}_{FSI}} J \boldsymbol{\sigma}_{f} \mathbf{F}^{-T} \boldsymbol{n} \cdot (\mathbf{v}_{s}-\mathbf{v}_{f}) \ d A_{\hat{x}}+ \int_{\hat{\Gamma}_{in}} p_{in}\boldsymbol{n} \cdot \mathbf{v}_{f} \ d A_{\hat{x}}\\
    &+\int_{\hat{\Gamma}_{out}} p_{out}\boldsymbol{n} \cdot \mathbf{v}_{f} \ d A_{\hat{x}}-\int _{ \hat{\Omega}_{s}}  \mathbf{P} : \nabla \mathbf{v}_{s}\ d \hat{x}\\
    &   - \int _{ \hat{\Omega}_{s}}\boldsymbol{b} _{s} \cdot \mathbf{v}_{s}\ d \hat{x}  -\int _{ \hat{\Omega}_{f}}J \boldsymbol{b} _{f}\cdot \mathbf{v}_{f}\ d \hat{x}=0 \\
    &\quad \text{for every } (\mathbf{v}_{f},\mathbf{v}_{s}) \in V_{f} \times V_{s}.\\
    \end{aligned}
\end{equation}
\item \textbf{FLUID INCOMPRESSIBILITY EQUATION:}
    \begin{equation}
        - \int _{ \hat{\Omega}_{f}}\operatorname{div}(J \mathbf{F}^{-1} \mathbf{u}_{f}) q_{f}\ d \hat{x} =0 \quad \text{for every } q_{f} \in Q_{f}.
    \end{equation}
\item \textbf{EXTENSION EQUATION}:
    \begin{equation}
        \label{extensionequation}
        \int _{ \hat{\Omega}_{f}} \nabla \boldsymbol{d} _{f}: \nabla \boldsymbol{e} _{f} \ d \hat{x} =0 \quad \text{for every } \boldsymbol{e}_{f} \in E_{f}.
    \end{equation}
\end{itemize}
The previous system can be approached numerically using numerous strategies. Our numerical resolution was built upon a monolithic method, which consists in solving the fully coupled problem, without having to iterate between the two sub-modules (which is instead done with the segregated approach \cite{degrote}). The advantages of this kind of treatment are the achieved stability properties \cite{stabilitymonolithic}, and the possibility to solve the whole problem employing a global nonlinear solver such as Newton linearization.
The price to pay is in terms of computational cost, which is why we combined this approach with the use of a modeling reduction strategy (Section \ref{reducedordermodel}).
\section{Numerical Approximation}%
\label{sec:numerical_approximation}
In this section, we will present the ingredients needed to tackle numerically a generic bifurcation problem held by a nonlinear PDE.
The starting point will be represented by the so-called full order (or high-fidelity) discretization, which aims at providing a high accuracy approximation to the solution of the PDE. With this regard, we will present only the Galerkin-Finite Element (FE) method which has been used for our simulation. However, one could employ our framework and change the underline full order discretization strategy.

We will present the algorithm we used for the reconstruction of the bifurcation diagram. The chosen approach is branch-wise because it aims at reconstructing the different branches separately. It also requires coupling the numerical resolution with a Continuation method\cite{continuation}, which fosters convergence to solutions belonging to the branch to be approximated.

The mitigation of the computational cost required by the presented methodology will be tackled using a reduced order modelling technique.

\subsection{Full order approximation}%
\label{sub:full_order_approximation}
The Galerkin-FE method is a numerical discretization technique that allows solving various PDEs by resorting to their variational formulation.
Here, we used this method to solve the FSI problem (Section \ref{sub:variational_formulation}), but it will also assume a prominent role in the reduced order setting (Section \ref{reducedordermodel}).

We start by considering a finite dimensional subspace $V_{h} \subseteq V$ such that $\operatorname{dim}\left(V_{h}\right) \equiv N_{h}<+\infty$.
The main idea is to solve the weak variational formulation \eqref{variationalform} onto this generic subspace, that means: given  $\boldsymbol{\mu} \in \mathbb{P}$, find $X_{h}(\boldsymbol{\mu})$ in $\mathbb{X}$ s.t.\
\begin{equation}
    \label{variationalfor}
f(X_{h}(\boldsymbol{\mu}), Y_{h};  \boldsymbol{\mu})=0 \quad \forall\, Y_{h} \in \mathbb{X}_{h}.
\end{equation}
We will refer to the above problem as Discrete Weak Formulation ($\mathcal{DWF}$).

Since $\mathbb{X}_{h}$ has finite dimension, i.e.\ $\operatorname{dim}\left(\mathbb{X}_{h}\right) \equiv N_{h}<+\infty,$ then $\mathbb{X}_{h}=\operatorname{span}\{\varphi_{1}, \varphi_{2}, \ldots,$ $\varphi_{N_{h}}\}$, i.e.\ every element of $\mathbb{X}_{h}$ can be written as a linear combination of the basis functions $\{\varphi_{i}\}_{i=1}^{N_h}$.
Thus, we can express the solution $X_{h}(\boldsymbol{\mu})$  as $\mathbb{X}_{h} \ni X_{h}(\boldsymbol{\mu})=\sum_{j=1}^{N_{h}} X_{j}(\boldsymbol{\mu}) \varphi_{j}$.

Because of the former consideration, one can easily prove that solving ($\mathcal{DWF}$) is equivalent to finding the solution $\boldsymbol{X_{h}}(\boldsymbol{\mu})=\{X_{j}(\boldsymbol{\mu})\}_{j=1}^{N_{h}} \in \mathbb{R}^{N_{h}}$ of the following nonlinear system:
\begin{equation}
    \label{nonlinsys}
f\left(\sum_{j=1}^{N_{h}} X_{j}(\boldsymbol{\mu}) \varphi_{j}, \varphi_{i}\right)=0, \quad \forall \ i=1\ldots N_{h}.
\end{equation}
The nonlinearity in the model forces us to consider a non-linear solver for its approximation, thus we relied on the Newton method.

First, we introduce a parametrized linear variational form by making use of the Fr\'echet partial derivative of $F$ at $Z \in \mathbb{X}$ w.r.t.\ $X$ as
\begin{equation}
    \mathrm{d} f[Z](X, Y ; \boldsymbol{\mu})\doteq\left\langle D_{X} F(Z ; \boldsymbol{\mu}) X, Y\right\rangle \quad \forall\ X, Y \in \mathbb{X}.
\end{equation}
Then, one suitably chooses an initial guess $X_{0} \in \mathbb{X}_{h}$, and proceed with a sequence of iterations composed by:
\begin{enumerate}
    \item Solving the finite dimensional linear variational problem:
        \begin{equation*}
            \mathrm{d} f\left[{X}_{h}^{k}(\boldsymbol{\mu})\right]\left(\delta X_{h}, Y_{h} ; \boldsymbol{\mu}\right)=-f\left({X}_{h}^{k}(\boldsymbol{\mu}), Y_{h} ; \boldsymbol{\mu}\right), \quad \forall\ Y_{h} \in \mathbb{X}_{h},
        \end{equation*}
   that can be cast in the following matrix formulation:
    $$
\mathbf{J}_{h}\left(\mathbf{X}_{h}^{k} ; \boldsymbol{\mu}\right) \delta \mathbf{X}=-\mathsf{F}_{h}\left(\mathbf{X}_{h}^{k} ; \boldsymbol{\mu}\right),
$$
being $\mathbf{J}_{h}$ the Jacobian matrix and $\mathsf{F}_h$ the high-order residual vector:
$$
\begin{aligned}
    (\mathbf{J}_{h})_{i j}:=\mathrm{d} f\left[{X}_{h}^{k}(\boldsymbol{\mu})\right]\left(\varphi_{j}, \varphi_{i}\right)  \quad \forall i,j&=1, \ldots, N_{h}, \quad \mathbf{J}_{h} \in \mathbb{R}^{N_{h} \times N_{h}}, \\
(\mathsf{F}_{h})_i:=f\left({X}_{h}^{k}(\boldsymbol{\mu}), \varphi_{i} ; \boldsymbol{\mu}\right)  \quad \forall i&=1, \ldots, N_{h}, \quad \mathsf{F}_{h} \in \mathbb{R}^{N_{h}}.\end{aligned}
$$ \item Update the approximated solution:
        \begin{equation}
X_{h}^{k+1}(\boldsymbol{\mu})=X_{h}^{k}(\boldsymbol{\mu})+\delta X_{h} .
\end{equation}
\end{enumerate}
As stopping criterion for the former algorithm one can use a control on the high-fidelity residual and/or on the norm of the variation $\delta X_{h}$, combined with a check on the maximum number of iterations.

We now want to use the former numerical resolution to perform a branch-wise reconstruction of the bifurcation phenomena \cite{Pichi}. This means following a single branch while modifying the value of the parameter we are dealing with. However, after the bifurcation had occurred, we lose the uniqueness of the solution, and therefore this task is not trivial. This is why we will implement a \textit{continuation method}, that is a procedure that allows to generate a sequence of solutions associated with the same branch \cite{continuation}.
We consider a discrete version of the parameter space $ \mathbb{P} _{h} \subset \mathbb{P} $, and let us imagine starting from a certain value of the  parameter $ \boldsymbol{\mu}_j\in \mathbb{P}_{h} $, for which the corresponding full order solution ${X}_{h}\left(\boldsymbol{\mu}_{j}\right)$ is known.
The idea is to get the solution ${X}_{h}\left(\boldsymbol{\mu}_{j+1}\right)$ for the next value of the parameter $\boldsymbol{\mu}_{j+1} \in \mathbb{P}_{h}$, starting from a suitable initial guess given by ($X_{h}(\boldsymbol{\mu}_{j})$, $\boldsymbol{\mu}_{j}$) .
In the following, we will be using the simplest version of this methodology which consists in choosing directly ($ X_{h}(\boldsymbol{\mu}_{j})$, $\boldsymbol{\mu}_{j}+ \Delta \boldsymbol{\mu}$) as initial guess. However, there are more complex strategies which can be exploited to automatically discover new branches \cite{classic_deflation, pintore2019efficient}.
We remark that a lot of attention must be paid to the choice of the discretization step; if it is too large, we risk not to capture the bifurcating behaviour, while a too-small step causes a possible waste of computational resources, especially in regions of the parameter space far from the bifurcation points.
The former considerations can be summarized in the following pseudo-algorithm for the reconstruction of a branch (see \cite{Pichi} for more details):
\begin{algorithm}[H]
\caption{A pseudo-code for the reconstruction of a branch}
\begin{algorithmic}[1]
\State{$\mathbf{X}_{h,0}=\mathbf{X}_{h,guess}$}\Comment{Initial guess}
\For{$\boldsymbol{\mu}_j \in \mathbb{P}_{h}$}\Comment{Continuation loop}
\State{{$\mathbf{X}_{h,j}^{0}=\mathbf{X}_{h,j-1}$}} \Comment{{Continuation guess}}
	\While{$||\mathsf{F}_{h}\left(\mathbf{X}_{h,j}^{k} ; \boldsymbol{\mu}_{j}\right)||_{\mathbb{R}^{N_{h}}} > \epsilon$}\Comment{Newton method}
    \State{$\mathbf{J}_{h}\left(\mathbf{X}_{h,j}^{k} ; \boldsymbol{\mu}_{j}\right) \delta \mathbf{X}_h=-\mathsf{F}_{h}\left(\mathbf{X}_{h,j}^{k} ; \boldsymbol{\mu}_{j}\right)$}
    \State{$\mathbf{X}_{h,j}^{k+1} = \mathbf{X}_{h,j}^{k}+ \delta \mathbf{X}_h$}
        \EndWhile
\EndFor
\end{algorithmic}
\label{algorithm-fom}
\end{algorithm}

\subsection{Reduced order model}
\label{reducedordermodel}
When dealing with a class of parametrized problem in the form:
\begin{equation}
    \label{parmStrong}
    F(X(\boldsymbol{\mu}); \boldsymbol{\mu})=0 \quad \text{with} \quad \boldsymbol{\mu} \in \mathbb{P} \subset \mathbb{R}^{p},
\end{equation}
we may be interested in looking for the solution for all the possible instances of the parameters' value, that is finding $X=X( \boldsymbol{\mu}) \in \mathbb{X}$  for all $\boldsymbol{\mu} \in \mathbb{P}$.
Basically, we would like to approximate the solution manifold $ \mathcal{S}$ (Equation \eqref{solutionmanif}).
So far, we have seen how to tackle the problem within a high-order approximation framework (Section \ref{sub:full_order_approximation}). However, given that the full order models present a considerable computational difficulty (which scales proportionally to the number of degrees of freedom $N_{h}$), various Reduced Order Models (ROMs) \cite{benner2017model, morepas2017} have been proposed to tackle the mitigation of the computational cost.

In the present work we have decided to use a particular ROM known as the \textit{Reduced Basis} (RB) method \cite{hesthaven2015certified,patera07:book,QuarteroniManzoniNegri2015}.
In particular, the RB methodology is divided into two stages:
\begin{itemize}
    \item \textbf{Offline phase}: high-fidelity numerical simulations for different parameter values are performed. These values must be representative, and several methods specialize in how that selection is made. This phase represents the most expensive one from the computational point of view, but it has to be done only once. The information obtained must be subsequently processed to be transmitted to the next phase.
    \item \textbf{Online phase} The precomputed quantities of the offline phase are used to solve with reduced computational complexity new instances of the parameterized problem. Ideally, the new resolutions' complexity should not depend on $N_{h}$ and should require resources that can be managed by low-power devices such as smartphones.
\end{itemize}

The previous strategy involves building a suitable approximation of the solution manifold $ \mathcal {S} $ by exploiting the information collected during the offline phase.
In particular, the RB method consists in constructing this approximation through a finite dimensional space $ \mathbb {X} _ {rb} $ generated by a suitable basis\footnote{Here, we will denote with the subscript ${}_{rb}$ all the reduced basis quantities.}.
In the next section, we will introduce the optimality condition satisfied by $ \mathbb{X} _{rb} $, together with a strategy for its computation.

\subsubsection{POD-Galerkin method}
\label{podgalerkin}

As discussed before, a crucial point when applying the RB method is the construction of a proper basis during the offline phase.
In the present work, we use the Proper Orthogonal Decomposition (POD) \cite{POD}, which can be used to compress information in different settings by extracting a low dimensional representation of a given dataset.

One begins by introducing a discretization in the parameter domain, namely $\mathbb{P}_{train}\subset \mathbb{P}$, whose cardinality will be denoted by $N_{train}$.
If we solve problem \eqref{parmStrong} for a sufficiently rich discretized parameter space $\mathbb{P}_{train}$, we obtain an approximation of the high-order solution manifold:

\begin{equation}
    \mathcal{S} _{\delta}=\{X_{h}(\boldsymbol{\mu})\ |\ \boldsymbol{\mu} \in \mathbb{P}\}\subset \mathbb{X}_{h},
\end{equation}

given by the vector space
\begin{equation}
    \mathbb{X}_{N_{train}}=\operatorname{span}\{X_{h}(\boldsymbol{\mu})\ |\ \boldsymbol{\mu} \in \mathbb{P}_{train}\},
\end{equation}
identified by the span of the so called snapshots (i.e.\ truth solutions w.r.t.\ $\boldsymbol{\mu}$).
Then, the idea is to find the space $\mathbb{X}_{rb}$ that minimizes the following quantity
\begin{equation}
\sqrt{\frac{1}{N_{\text {train }}} \sum_{\boldsymbol{\mu} \in \mathbb{P}_{train}} \inf _{X_{rb} \in \mathbb{X}_{rb}}\left\|X_{h}(\boldsymbol{\mu})-X_{rb}\right\|_{\mathbb{X}}^{2}},
\end{equation}
among all the subspaces of $\mathbb{X}_{N_{train}}$ with fixed dimension $N_{rb}$.

From a discrete point of view, the former procedure is equivalent to carrying out a Single Value Decomposition (SVD) on the snapshot matrix
\begin{equation}
\boldsymbol{S}=\left[X_{h}\left(\boldsymbol{\mu}_{1}\right), \ldots, X_{h}\left(\boldsymbol{\mu}_{N_{\text {train}}}\right)\right] \in \mathbb{R}^{N_{h} \times N_{train}} .
\end{equation}
By using the first $N_{rb}$ left singular vectors as expansion coefficients w.r.t. $X_h(\mu_j)$, one obtains the set of RB functions $\{\Psi_{k}\}_{k=1} ^{N_{rb}}$ which span $\mathbb{X}_{rb}$.

The procedure in question must be carried out independently for every field of the unknown variable.
Once obtained the reduced basis space $\mathbb{X}_{rb}$, we perform a Galerkin projection onto the subspace itself.
This strategy is equivalent to the one presented in Section \ref{sub:full_order_approximation}, the only difference lays in the basis that we are using in this setting, which will be the one obtained via POD. Thus, the reduced problem reads: given $\boldsymbol{\mu} \in \mathbb{P}$, find $X_{rb}(\boldsymbol{\mu}) \in \mathbb{X}_{rb}$ s.t.
\begin{equation}
\left\langle F(X_{rb}(\boldsymbol{\mu}) ; \boldsymbol{\mu}), Y_{rb}\right\rangle \doteq f(X_{rb}(\boldsymbol{\mu}), Y_{rb}; \boldsymbol{\mu})=0 \quad \forall\ Y_{rb} \in \mathbb{X}_{rb} .
\end{equation}

The latter reduces, similarly to the full order model, to a system of $N_{rb}$ nonlinear equations, which can be solved by means of the Newton method, as presented in Section \ref{sub:full_order_approximation}. In particular, one can prove that the $k$-th  iteration of the Newton method can be cast in matrix form as
\begin{equation}
\mathbf{V}^{T} \mathbf{J}_{h}\left(\mathbf{V} \mathbf{X}_{rb}^{k}(\boldsymbol{\mu}) ; \boldsymbol{\mu}\right) \mathbf{V} \delta \mathbf{X}_{rb}=-\mathbf{V}^{T} \mathsf{F}_{h}\left(\mathbf{V} \mathbf{X}_{rb}^{k}(\boldsymbol{\mu}) ; \boldsymbol{\mu}\right),
\end{equation}
where $\mathbf{V} \in \mathbb{R} ^{N_{h}\times N_{rb}}$ is the matrix whose columns denote the coefficients of the reduced basis functions $\{\Psi_{k}\}_{k=1} ^{N_{rb}}$ in terms of the full order basis functions $\{\varphi_{i}\}_{i=1} ^{N_{h}}$.
With regard to the branch-wise reconstruction, one can obtain a reduced order version of Algorithm \ref{algorithm-fom}, by substituting the full order quantities with their reduced order counterpart:

\begin{algorithm}[H]
\caption{A pseudo-code for the reduced order reconstruction of a branch}
\begin{algorithmic}[1]
\State{$\mathbf{X}_{rb,0}=\mathbf{X}_{rb,guess}$}\Comment{Initial guess}
\For{$\boldsymbol{\mu}_j \in \mathbb{P}_{h}$}\Comment{Continuation loop}
\State{{$\mathbf{X}_{rb,j}^{0}=\mathbf{X}_{rb,j-1}$}} \Comment{{Continuation guess}}
	\While{$||\mathsf{F}_{rb}\left(\mathbf{X}_{rb,j}^{k} ; \boldsymbol{\mu}_{j}\right)||_{\mathbf{R}^{N_{rb}}} > \epsilon$}\Comment{Newton method}
    \State{$\mathbf{V}^{T} \mathbf{J}_{h}\left(\mathbf{V} \mathbf{X}_{rb,j}^{k} ; \boldsymbol{\mu}_{j}\right) \mathbf{V} \boldsymbol{\delta} \mathbf{X}_{rb,j}=-\mathbf{V}^{T} \mathsf{F}_{h}\left(\mathbf{V} \mathbf{X}_{rb,j}^{k} ; \boldsymbol{\mu}_{j}\right)$}
    \State{$\mathbf{X}_{rb,j}^{k+1}= \mathbf{X}_{rb,j}^{k}+ \delta \mathbf{X}_{rb}$}
        \EndWhile
\EndFor
\end{algorithmic}
\label{algorithm-reduced}
\end{algorithm}

\section{Numerical results}
\label{sec:numerical}
In this section, we will present the results of our simulations for the fluid-structure interaction problem. In particular, we will compare the multi-physics case and the rigid structure one, where a non-deformable rigid body replaces the elastic structure.
We will also present the results for the reduction strategy implemented through the RB method, whose accuracy has been tested on never seen parametric instances of the branch we aim at approximating.
To investigate the effect of large scale deformations, we will use the Saint Venant-Kirchoff nonlinear model for the solid constitutive relation.

Finally, all the proposed models will be compared, obtaining important insights on how the bifurcation point is modified by the introduction of the structure and the constitutive relation chosen to model it.

\subsection {Preliminaries}
\begin {figure}[h]
    \centering
    \includegraphics [width = 1 \linewidth] {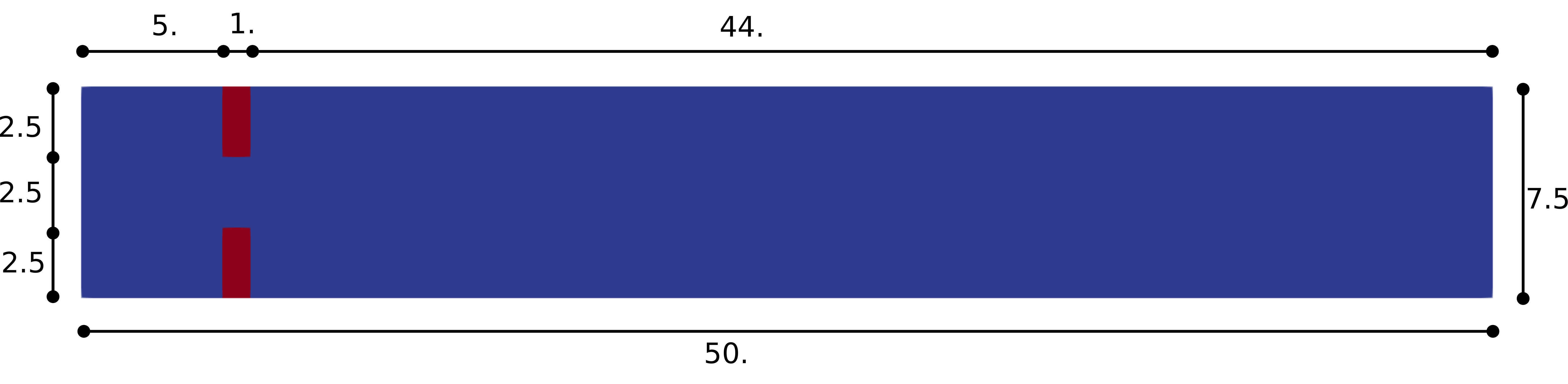}
    \caption {Computational domain which represents a straight channel with a sudden expansion.}%
    \label {fig: domainRes}
\end {figure}

Let us start by noting that we have used \textit{centimeter-gram-second (CGS)} unit system for all the equations.
Consequently, the various quantities have the following units of measurement:
\begin {table} [H]
\centering
\begin {tabular} {ll}
    \toprule
    \textbf {Physical quantity} & \textbf {Unit of measure}  \\
    \midrule
    lengths & cm \\
    $\sigma_{f},\mathbf{P}, p_{in}$ & barye (Ba)\\
    $\mathbf{u}$  & $\frac{cm}{s}$\\
    $\rho_{f},\rho_{s}$ & $\frac{g}{cm^3}$\\
    $ \mu$ & stokes (St) \\
    \bottomrule
\end {tabular}
\caption{Unit of measures of the main physical quantities involved in the FSI problem.}
\end{table}
The computational domain is the same for the FSI problem and the rigid structure one, and it is shown in Figure \ref{fig: domainRes}.
In the case of the rigid structure, the red region in Figure \ref {fig: domainRes} is removed from the computational domain and is used only to impose the boundary conditions. On the contrary, for the FSI problem, it represents the reference domain for the solid $ \hat {\Omega} _ {s} $ (see Figure \ref{fig:domain}), where we solve the elastic problem .
As already anticipated in Section \ref{sub:full_order_approximation}, as high-fidelity approximation we adopted a FE discretization on unstructured meshes composed by triangular elements (see Table \ref{table1}).

\begin {table} [h]
\centering
\begin {tabular} {llcc}
    \toprule
    \textbf {Problem} & \textbf {Grid type} & \textbf {Number of cells} & \textbf {Number of vertices} \\
    \midrule
    FSI & Unstructured with triangles & 36419 & 18584 \\
    Rigid structure & Unstructured with triangles & 36118 & 18460 \\
    \bottomrule
\end {tabular}
\caption {Description of the computational meshes used for the simulations.}
\label{table1}
\end{table}

The numerical resolution was carried out using the FeniCS library \cite{fenics}, with the additional use of multiphenics \cite{multiphenics}, which facilitates the definition of multi-physics problems, and of RBniCS \cite{rbnics}, which instead implements several reduced order modeling techniques for parametrized problems.
\subsection{Rigid structure}%
\label{sec:rigid_structure}
The rigid structure case involves solving the Navier-Stokes system in Eulerian coordinates. In fact, since the domain is fixed there is no need to use the ALE formulation, so that the system of equations we solve is:
\begin{equation}
    \label{nschap4}
\left\{\begin{array}{ll}
        -\operatorname{div} \sigma_{f}(\mathbf{u}_f,p) +\rho_{f}(\mathbf{u}_f \cdot \nabla) \mathbf{u}_f=0 & \text { in } \hat{\Omega}_f , \\
\operatorname{div} \mathbf{u}_f=0 & \text { in } \hat{\Omega }_f, \\
\sigma_{f}\mathbf{n}=-p_{in} I & \text { on } \hat{\Gamma}_{in}, \\
\sigma_{f}\mathbf{n}=0 & \text { on } \hat{\Gamma}_{out}, \\
\mathbf{u}_f=\mathbf{0} & \text { on } \Gamma_{0} ;
\end{array}\right. \quad \text{where} \quad
    \begin{dcases}
        \hat{\Gamma}_{in}={0} \times [0,7.5],\\
        \hat{\Gamma}_{out}={50} \times [0,7.5],\\
        \hat{\Gamma}_0=\hat{\Omega}_f \backslash \{ \hat{\Gamma}_{in} \cup \hat{\Gamma}_{out}\}.
    \end{dcases}
\end{equation}

The boundary conditions  on $ \hat{\Gamma}_ {in} $ and $ \hat{\Gamma}_ {out} $ are of Neumann type and concern the imposition of the inlet and outlet stress values. In particular, the condition on $ \hat{\Gamma}_{out} $ is stress-free being homogeneous.
Therefore, the flow is exclusively animated by the value of  $ p_ {in} $, which encodes the gradient of the stresses between the inlet and the outlet. This choice represents an alternative to a Dirichlet boundary condition which directly imposes the velocity value at the inlet. This last approach has been excluded to avoid lifting technique and make the ROM approximation easier.
For our simulations we set $ p_ {in} = 450$, and $\rho_{f}=1$, while the parameter adopted for the bifurcation analysis is the kinematic viscosity\footnote{Since we have written the equations in the $CGS$ system, and being $\rho_{f}=1$, the values of the kinematic viscosity and those for dynamic viscosity coincide (except for the dimensions).} of the fluid $\mu$.
However, since the channel's geometry is fixed, and the only physical parameter that varies is the viscosity, the dimensionless bifurcation parameter is the Reynolds number.
For what concerns the FE discretization, we used the Taylor-Hood ($\mathcal{P   }_2-\mathcal{P}_1$)  elements to avoid stability issues \cite{quarteroni2008numerical}.

\begin{figure}[h]
    \centering
    \includegraphics[width=0.75\linewidth]{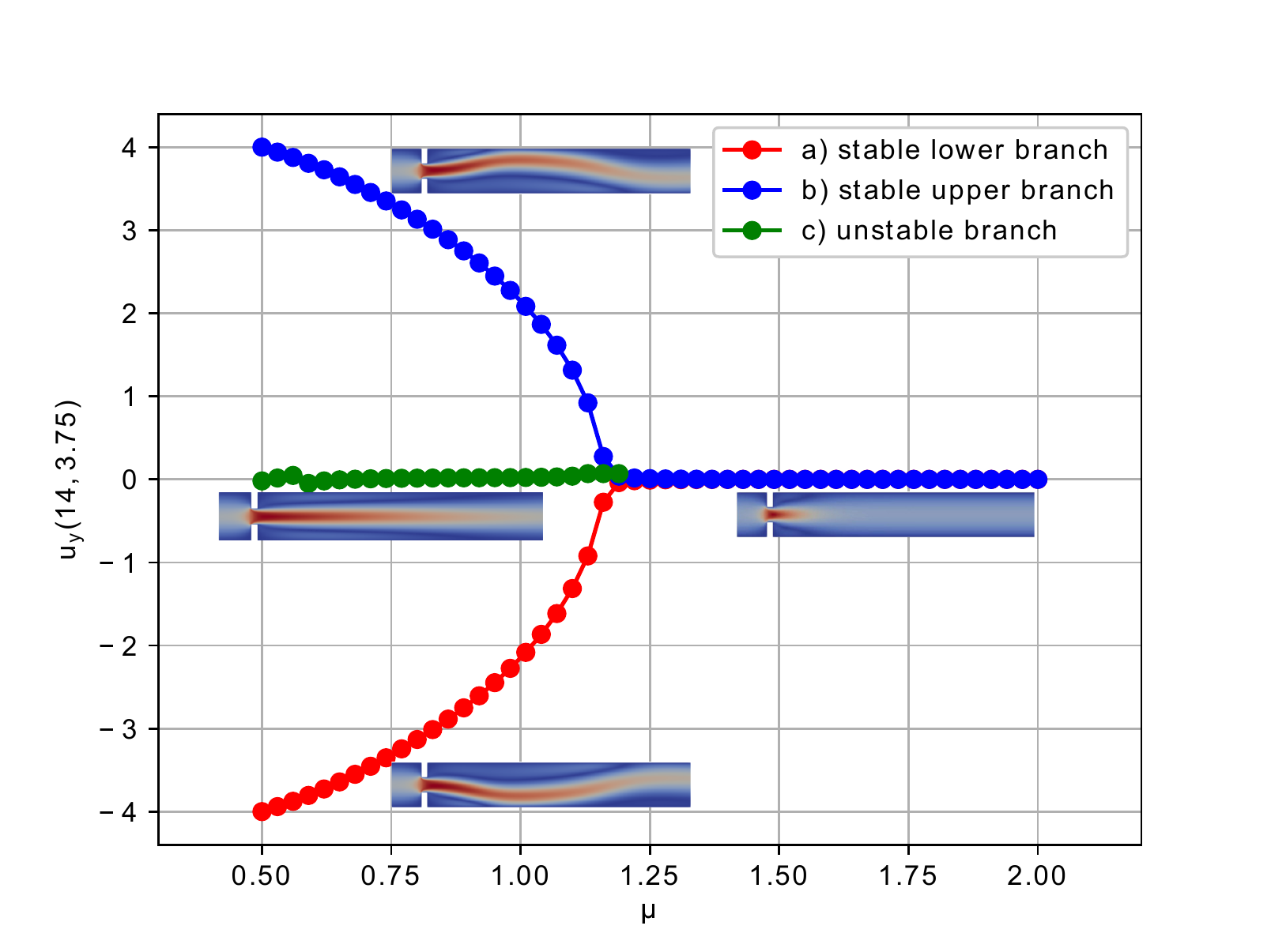}
    \caption{Bifurcation diagram for the Navier-Stokes system.}%
    \label{bifNSystem}
\end{figure}
We varied the kinematic viscosity $ \mu $ within the parametric interval $ \mathbb {P}_h = [0.5,2] $, composed by $51$ equispaced points (which is equivalent to adopt a step $\Delta \mu=0.03$).

During the study of the system, we expect a behavior similar to the one investigated in \cite {Pichi, coanda}, that is a \textit {supercritical pitchfork bifurcation}. This consists in having a stable symmetric solution above a critical value $ \mu ^ {*} $, while, below this value, two stable wall-hugging solutions and an unstable symmetric one coexist.
To provide the bifurcation diagram shown in Figure \ref{bifNSystem}, we have chosen as a representative scalar quantity the vertical velocity at a point of the channel symmetry axis, that is $ u_ {y} (\underline {\boldsymbol {x}}) $, with $ \underline {\boldsymbol {x}} = (14,3.75) $.
It should be noted that the procedure implemented to follow both the unstable symmetric branch and the two stable asymmetric branches is the same (see Algorithm \ref{algorithm-fom}); the only difference, within our approach, consists in the direction of exploration of the parameter space $\mathbb{P}_{h}$.
The two branches (a), (b) of Figure \ref{bifNSystem} refer to the asymmetric wall-hugging jets, respectively to the upper and lower walls. In order to reconstruct these two branches, we proceed by calculating the solution for $ \mu = 2.0 $ which is obtained starting from the trivial null guess, and then we decrease the viscosity through the $ \Delta \mu $ step. However, due to the unstable nature of branch (c), for its reconstruction it is necessary to start with $\mu=0.5$ and proceed by increasing instead of decreasing it.

\begin{figure}
\centering
\subfloat[$\mu=0.98$]{%
  \includegraphics[width=1\textwidth]{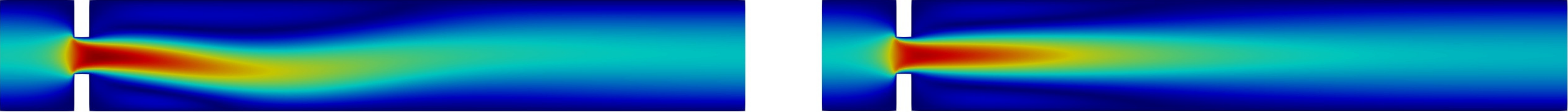}%
  }\par
\subfloat[$\mu=0.71$]{%
  \includegraphics[width=1\textwidth]{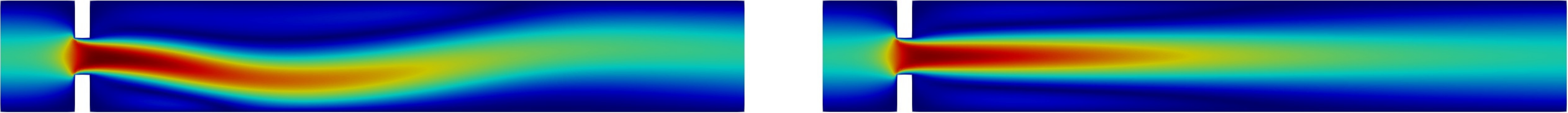}%
  }\par
\subfloat[$\mu=0.5$]{%
  \includegraphics[width=1\textwidth]{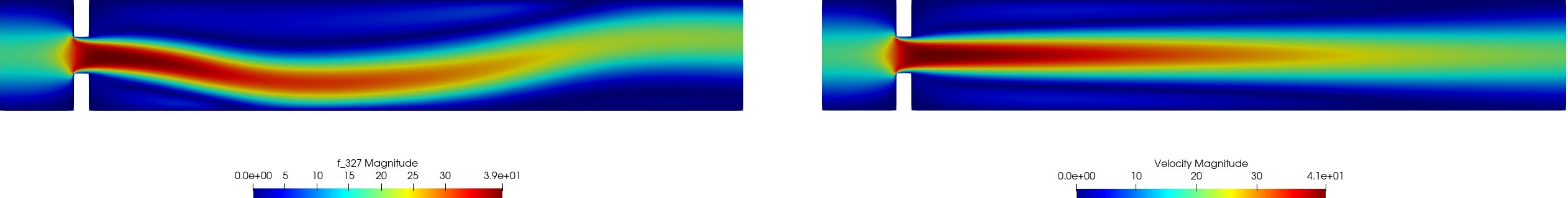}%
  }
\caption{Velocity magnitude for some significant values of the kinematic viscosity: the left solutions belong to the lower stable branch (a), the right ones to the unstable middle branch (c).}
\label{fig:rigidsnapvel}
\end{figure}

\begin{figure}
\centering
\subfloat[$\mu=0.98$]{%
  \includegraphics[width=1\textwidth]{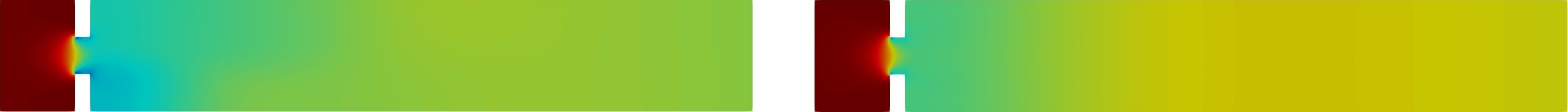}%
  }\par
\subfloat[$\mu=0.71$]{%
  \includegraphics[width=1\textwidth]{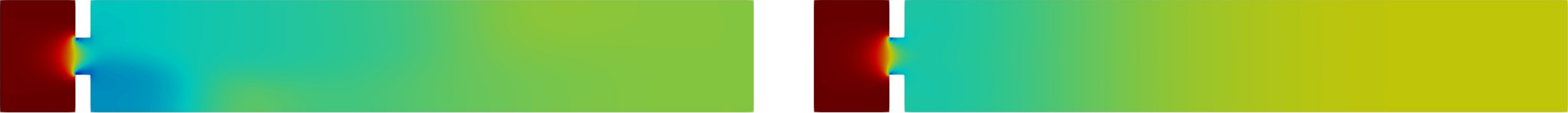}%
  }\par
\subfloat[$\mu=0.5$]{%
  \includegraphics[width=1\textwidth]{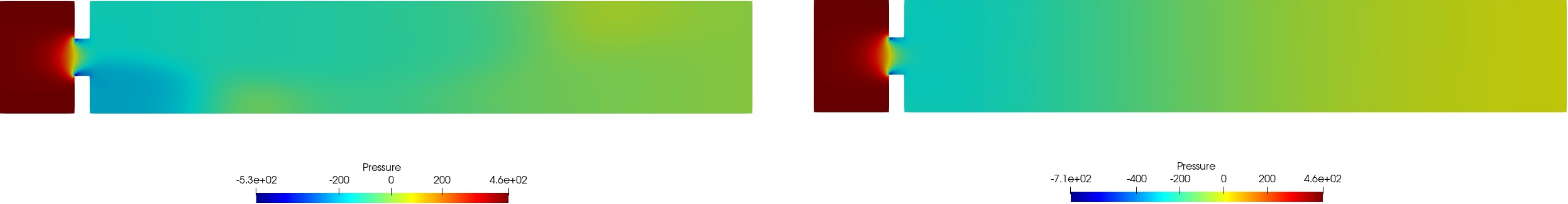}%
  }
\caption{Pressure field for some significant values of the kinematic viscosity: the left solutions belong to the lower stable branch (a), the right ones to the unstable middle branch (c).}
\label{fig:rigidsnappres}
\end{figure}

\begin{figure}[H]
    \centering
    \includegraphics[width=1\linewidth]{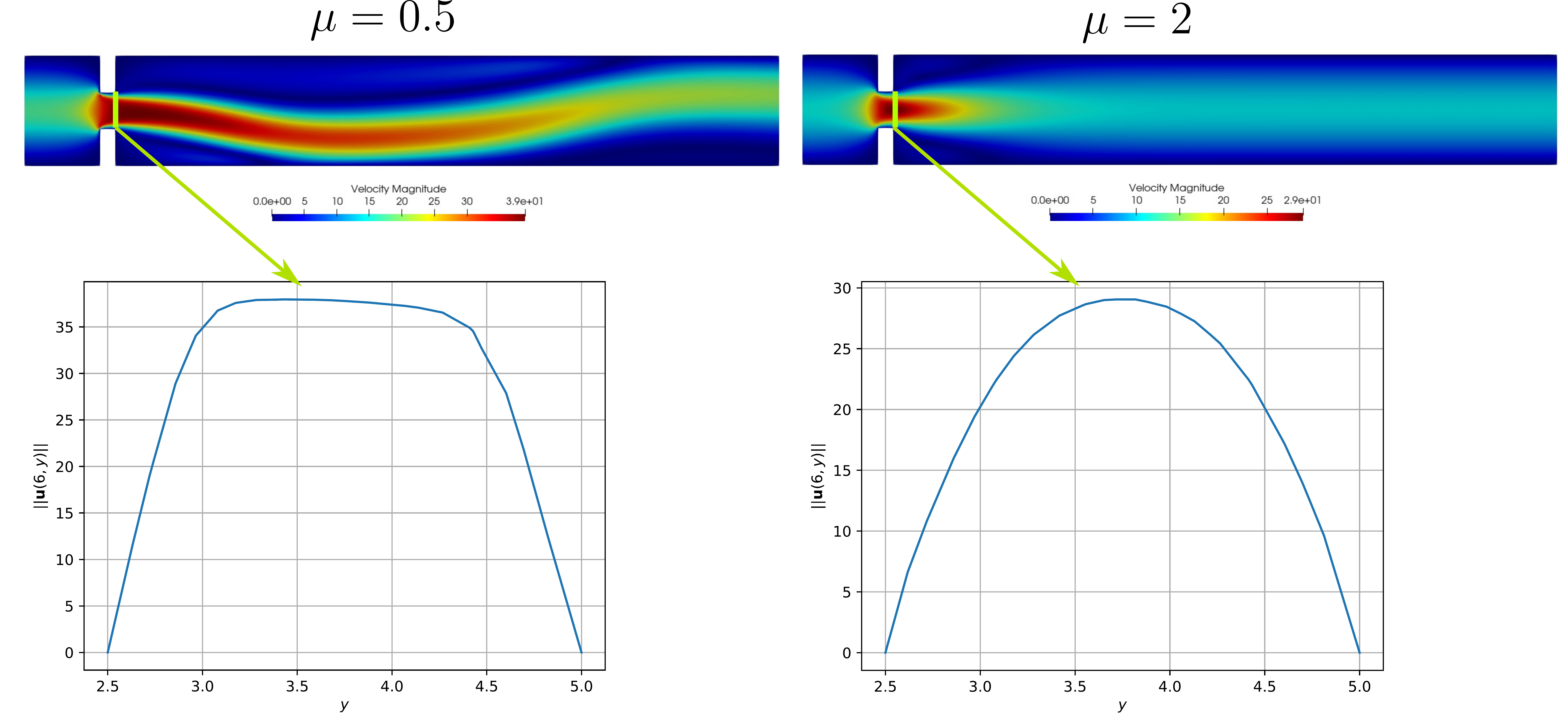}
    \caption{Velocity magnitude profiles at the expansion section $\{6\} \times [2.5,5]$ .}%
    \label{fig:velocity-inlet}
\end{figure}

As shown in Figure \ref{bifNSystem}, we were able to completely reconstruct the bifurcation diagram and estimate that, for this particular setting, the bifurcation occurs at $ \mu ^ {*} \approx 1.19 $. Figures \ref{fig:rigidsnapvel} and \ref{fig:rigidsnappres} show the velocity and pressure profiles belonging to the branches $(a) $ and $(c)$ for values of $ \mu< \mu ^{*} $.
It is important to note that although the inlet stress $ p_ {in} $ is constant, the velocity at the exit of the  expansion orifice varies according to the viscosity. This is highlighted in Figure \ref{fig:velocity-inlet}, which shows the velocity profile  for the extreme values of the parametric range, namely $\mu \in\{0.5,2.0\} $
\captionsetup[subfigure]{position=top}


\captionsetup[subfigure]{position=bottom}
%
\subsection{FSI problem with linear elasticity}%
\label{sec:fsi_problem_with_linear_elastisticy}
In this section, we introduce the results related to the fluid-structure interaction problem. We will start by replacing the rigid structure with a linear elastic one composed by two leaflets.  We recall that the FSI system we need to solve is

\begin{equation}
    \label{steadyfsilinear}
\left\{\begin{array}{ll}
        \rho_{f}J { \nabla}{\mathbf{u}}_{f}\mathbf{F}^{-1}{\mathbf{u}}_{f}   -{\operatorname{div}}(J {\boldsymbol{\sigma}}_{f} ({\mathbf{u}}_{f}, p_{f})\mathbf{F}^{-T}) =0 \quad & \text{in} \quad \hat{\Omega}_{f},\\
        {\operatorname{div}} (J \mathbf{F}^{-1}\mathbf{u}_{f})=0 \quad & \text{in} \quad \hat{\Omega} _{f}, \\
        {\Delta} \boldsymbol{d}_{f}=0  \quad  & \text{in} \quad \hat{\Omega}_{f},\\
        -{\operatorname{div}} \mathbf{P} ({\boldsymbol{d} _{s}}) =0 \quad & \text{in} \quad \hat{\Omega}_{s};
 \end{array}\right.
\end{equation}
which is completed by the constitutive relations,  boundary conditions and  coupling conditions introduced in Section \ref{strongformulationFSI}.
An important remark concerns the imposition of the coupling conditions,
In fact we decided to weakly impose the constraints $\mathbf{u}_{f}=0$ and $\boldsymbol{d}_{f}=\boldsymbol{d}_{s}$, by means of a \textit{Lagrange multipliers} \cite{YU20051} approach, which implies introducing two new unknown $\boldsymbol{\lambda}_{d},\boldsymbol{\lambda}_{u}\in H^{1}(\hat{ \Gamma}_{FSI};\mathbb{R}^{2})$.
The space discretization is obtained using second order Lagrange FE for $\mathbf{u}_{f},\boldsymbol{d}_{f}$ and $\boldsymbol{d}_{s}$, while first order Lagrange FE for $p,\boldsymbol{\lambda}_{u}$ and $\boldsymbol{\lambda}_{d}$.

The solid's constitutive relation is completely determined by the value of the Lamè constants. These two parameters have the same dimension as stress, and for the results of this section, they assume the following values:

\captionsetup[subfigure]{position=top}

\begin{figure}[h]
    \centering
    \includegraphics[width=0.5\linewidth]{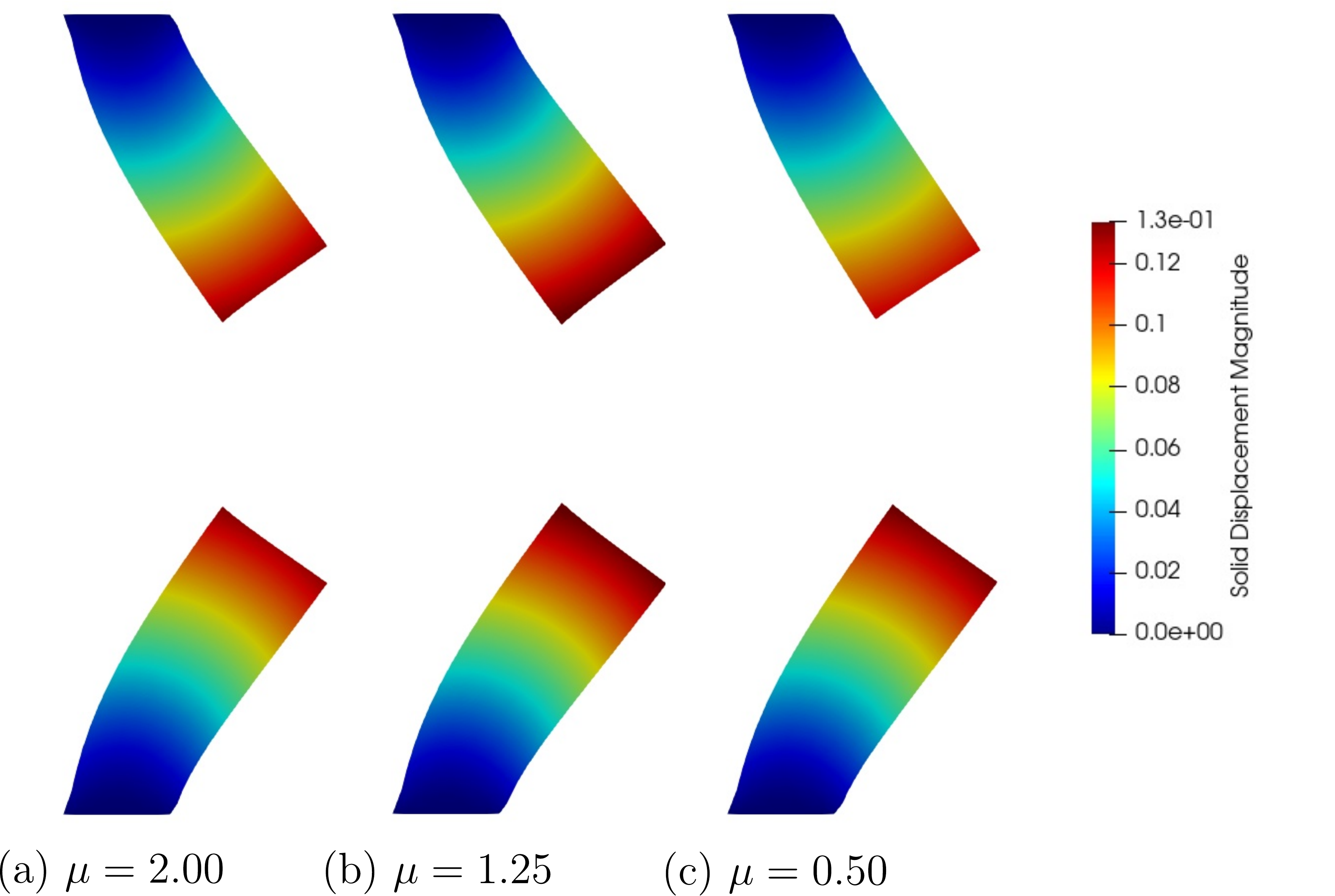}
    \caption{Solid's deformation snapshots: the solutions are given in a spatial frame through the $ \mathcal{A} _{f}$ map. }%
    \label{fig:snapshots-fsi-disp}
\end{figure}

\begin {table} [H]
\centering
\begin {tabular} {lc}
    \toprule
    \textbf {Physical quantity} & \textbf {Value}  \\
    \midrule
    $\lambda_{s}$ (Lamé's first parameter) & $8 \cdot 10^{5}$ Ba\\
    $\mu_{a}$ (Shear modulus) & $10^{5}$ Ba \\
    \bottomrule
\end {tabular}
\caption{}
\end{table}
\begin{figure}[H]
\centering
\subfloat[$\mu=2.00$]{%
  \includegraphics[width=1\textwidth]{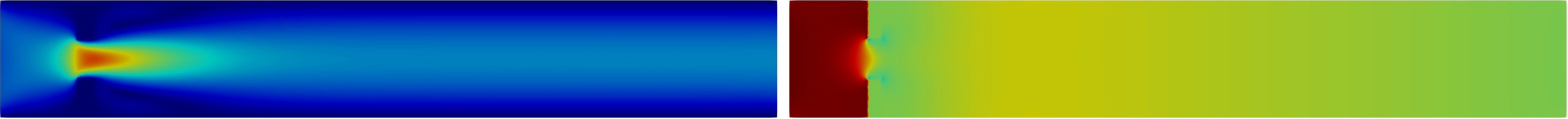}%
  }\par
\subfloat[$\mu=1.70$]{%
  \includegraphics[width=1\textwidth]{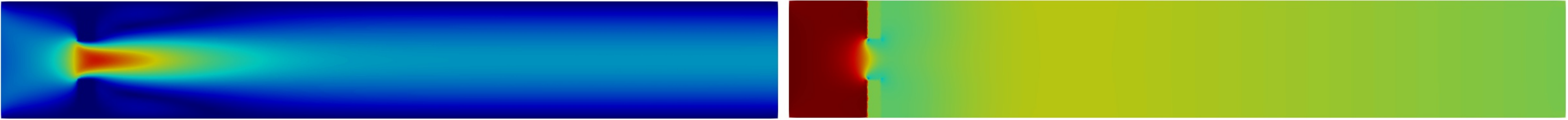}%
  }\par
\subfloat[$\mu=1.40$]{%
  \includegraphics[width=1\textwidth]{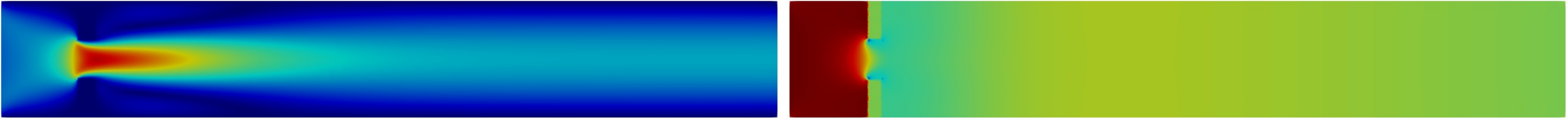}%
  }\par
\subfloat[$\mu=1.10$]{%
  \includegraphics[width=1\textwidth]{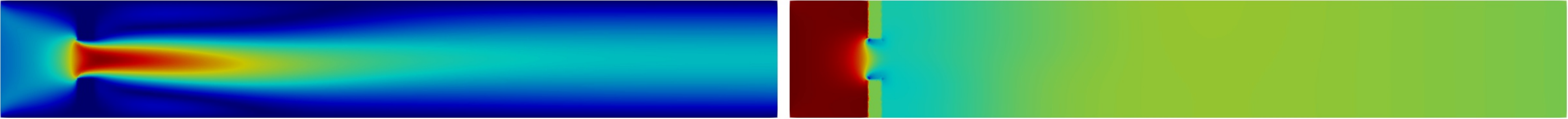}%
  }\par
\subfloat[$\mu=0.80$]{%
  \includegraphics[width=1\textwidth]{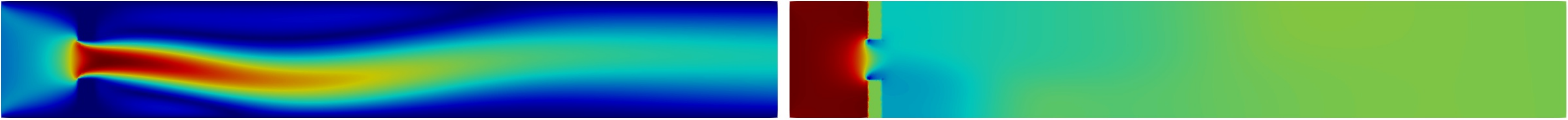}%
  }\par
\subfloat[$\mu=0.50$]{%
  \includegraphics[width=1\textwidth]{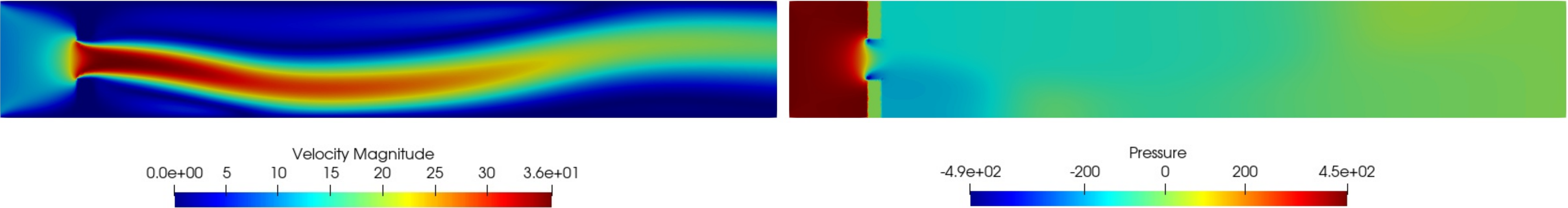}%
  }
  \caption{Velocity (left) and pressure (right) snapshots for different values of the kinematic viscosity. We can clearly notice the transversal gradient of pressure downstream of the expansion.}
  \label{fig:snapshots-fsi}
\end{figure}
\captionsetup[subfigure]{position=bottom}
The setting for the bifurcation analysis is the same as in Section \ref{sec:rigid_structure}, that is we consider the same values for the inlet and outlet stresses and we adopt the discretization $\mathbb{P}_{h}$, composed by $51$ equispaced points in the interval $[0.5,2.0]$. However, given the more complex nature of the problem, we have also recovered a RB approximation of the FSI bifurcating phenomena using the approach presented in Section \ref{reducedordermodel}.
First of all the offline phase was perfomerd by computing a FE solution to system \eqref{steadyfsilinear} for each $\mu_{j}$ in $\mathbb{P}_{h}$.
Figures \ref{fig:snapshots-fsi} and \ref{fig:snapshots-fsi-disp} show the high-fidelity solutions for solid's displacement, velocity, and  pressure obtained during this phase.

\begin{figure}[h]
    \centering
    \includegraphics[width=0.9\linewidth]{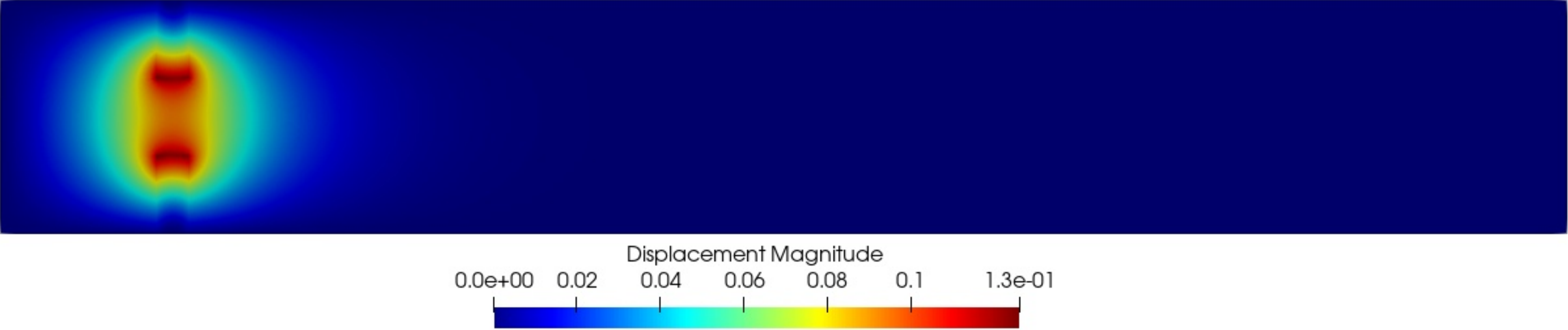}
    \caption{Global displacement $(\boldsymbol{d}_{f},\boldsymbol{d}_{s}$) for $\mu=2.0$  }%
    \label{fig:global-displacement}
\end{figure}
It should be noted that each solid displacement field is associated with a fluid displacement field which represents its harmonic extension. Indeed, if we consider the solution for $ \mu = 2.0 $ and plot $ \boldsymbol {d} _ {f} $ and $ \boldsymbol {d} _ {s} $ together we obtain the result shown in Figure \ref{fig:global-displacement}.
The structure is clearly that of a Poisson problem, and in fact $ \boldsymbol {d} _ {s} $, defined on $ \hat {\Omega} _ {s} $, spreads from the interface $ \hat {\Gamma} _ {FSI} $ in all $ \hat {\Omega} _ {f} $.

\begin{figure}[h]%
    \subfloat[\label{fig:bifurcation-fsi}Lower branch of the bifurcation diagram for the FSI problem.]{{\includegraphics[width=0.52\linewidth]{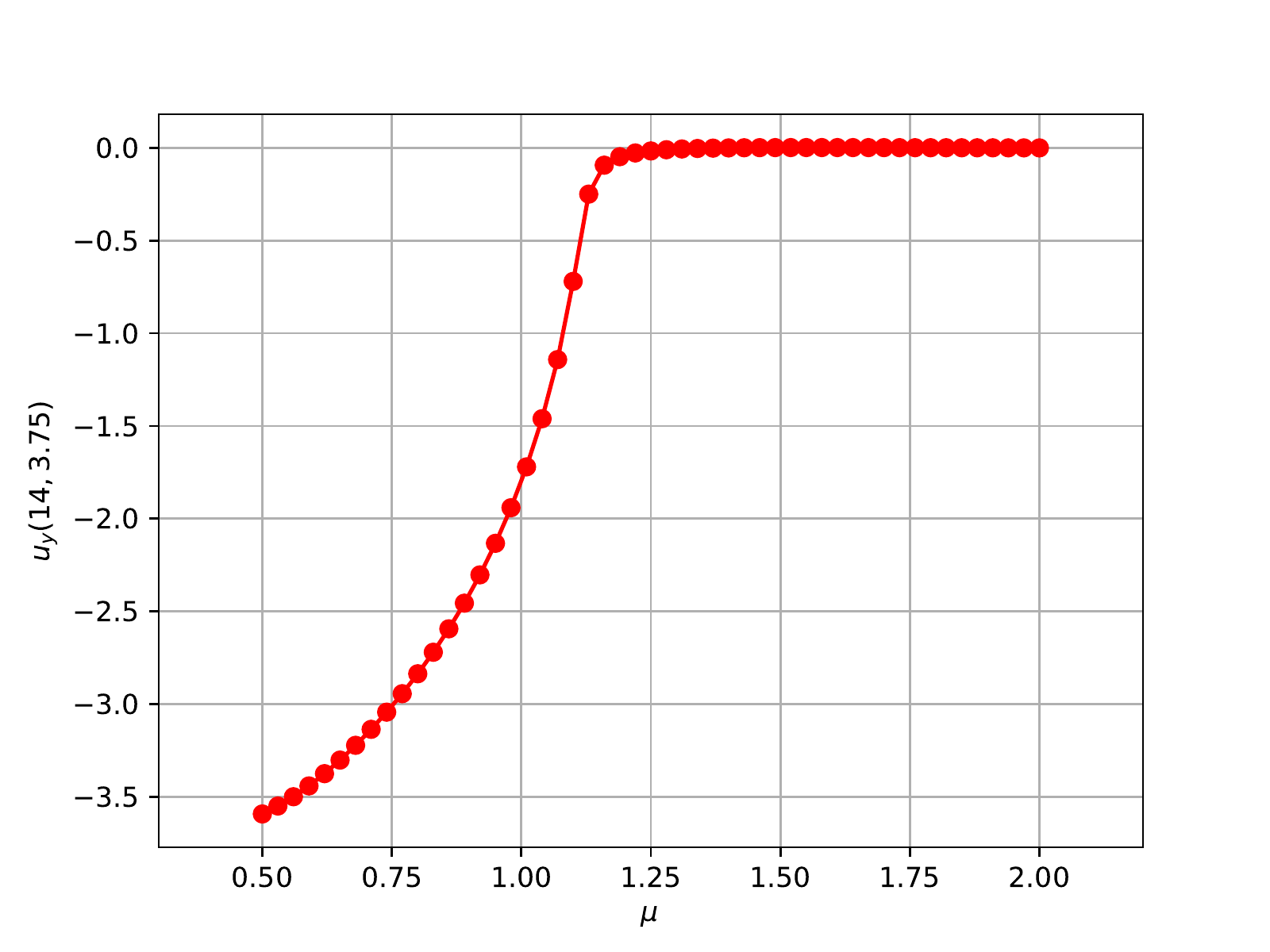} }}%
    \subfloat[\label{fig:max-displacement}Maximum magnitude of the solid's displacement of the upper and lower leaflets.]{{\includegraphics[width=0.52\linewidth]{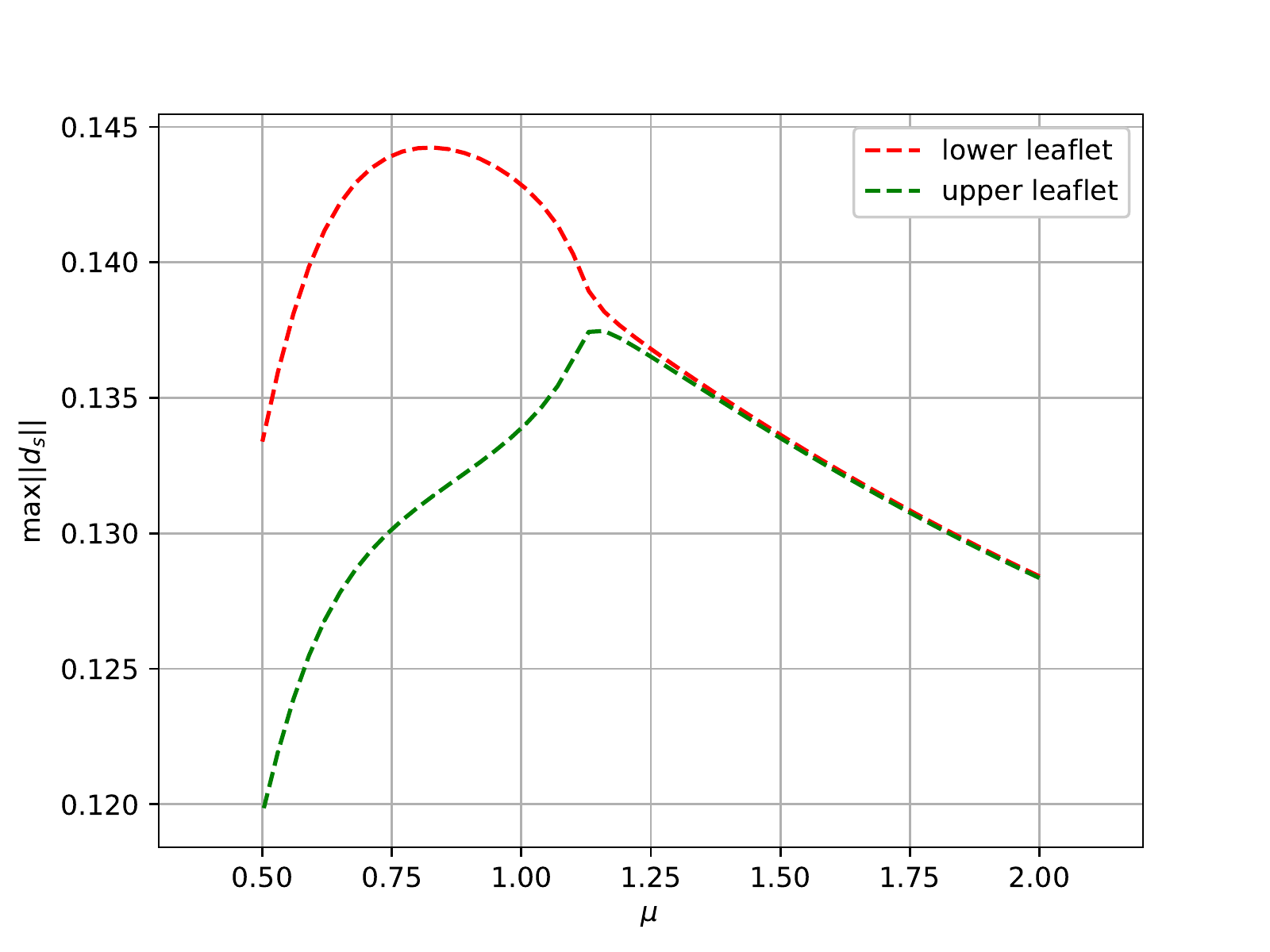} }}%
    \caption{}%
\end{figure}

We can now proceed as done in Section \ref {sec:rigid_structure} and retrieve a bifurcation diagram for the fluid-structure interaction problem. The first approach relies on considering the fluid phase and taking again as characteristic scalar output the vertical velocity $ u_ {y} (\underline {\boldsymbol {x}}) $, with $ \underline {\boldsymbol {x}} = ( 14,3.75) $.
If we limit ourselves to the lower stable branch, the result is shown in Figure \ref{fig:bifurcation-fsi}.

At this point, we want to investigate how the structure undergoes the bifurcation. This question arises spontaneously given that the symmetry breaking bifurcation is driven by the Navier-Stokes nonlinearity.
The first attempt consists in analyzing the behavior of the maximum displacement of the upper and lower leaflets. In fact, as long as the fluid flow is symmetric, we expect that the upper and lower displacements exhibit the same trend. However, when the flow becomes asymmetric, we expect the solid deformation field to become asymetric as well. Moreover the lower leaflet undergoes greater deformations, as the jet of fluid directed downwards creates a depression in the lower part of the channel.
Therefore, if we track $\max\|\boldsymbol{d}_{s}\|$ as a function of $\mu$ we obtain Figure\ref{fig:max-displacement}.

\begin{figure}[h]%
    \subfloat[\label{fig:deltap}Pressure drop through the upper and lower leaflets.]{{\includegraphics[width=0.5\linewidth]{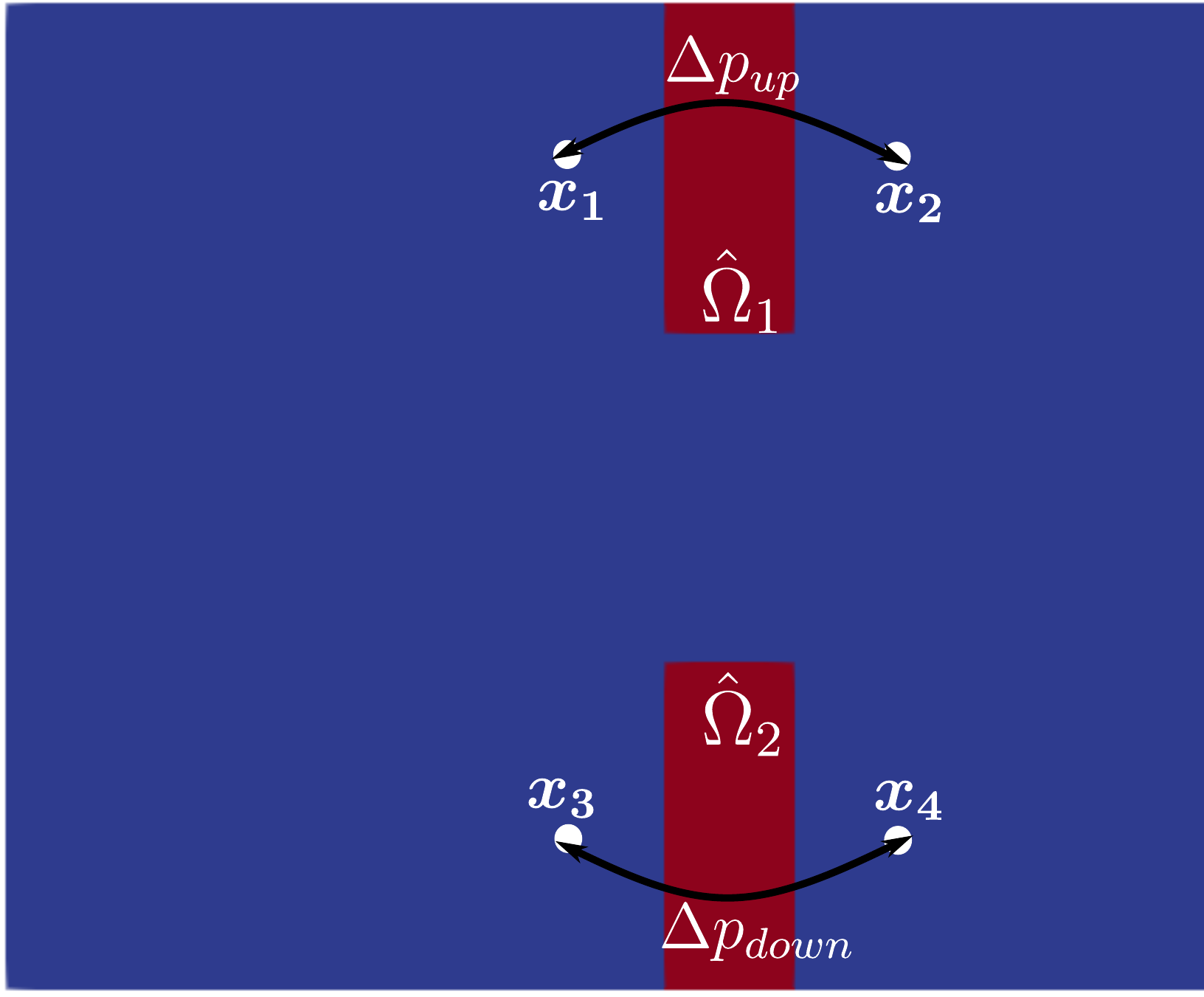} }}%
    \subfloat[\label{fig:deltaptrend}Pressure drop measured by $\Delta p_{up}$ and $\Delta p _{down}$.  ]{{\includegraphics[width=0.5\linewidth]{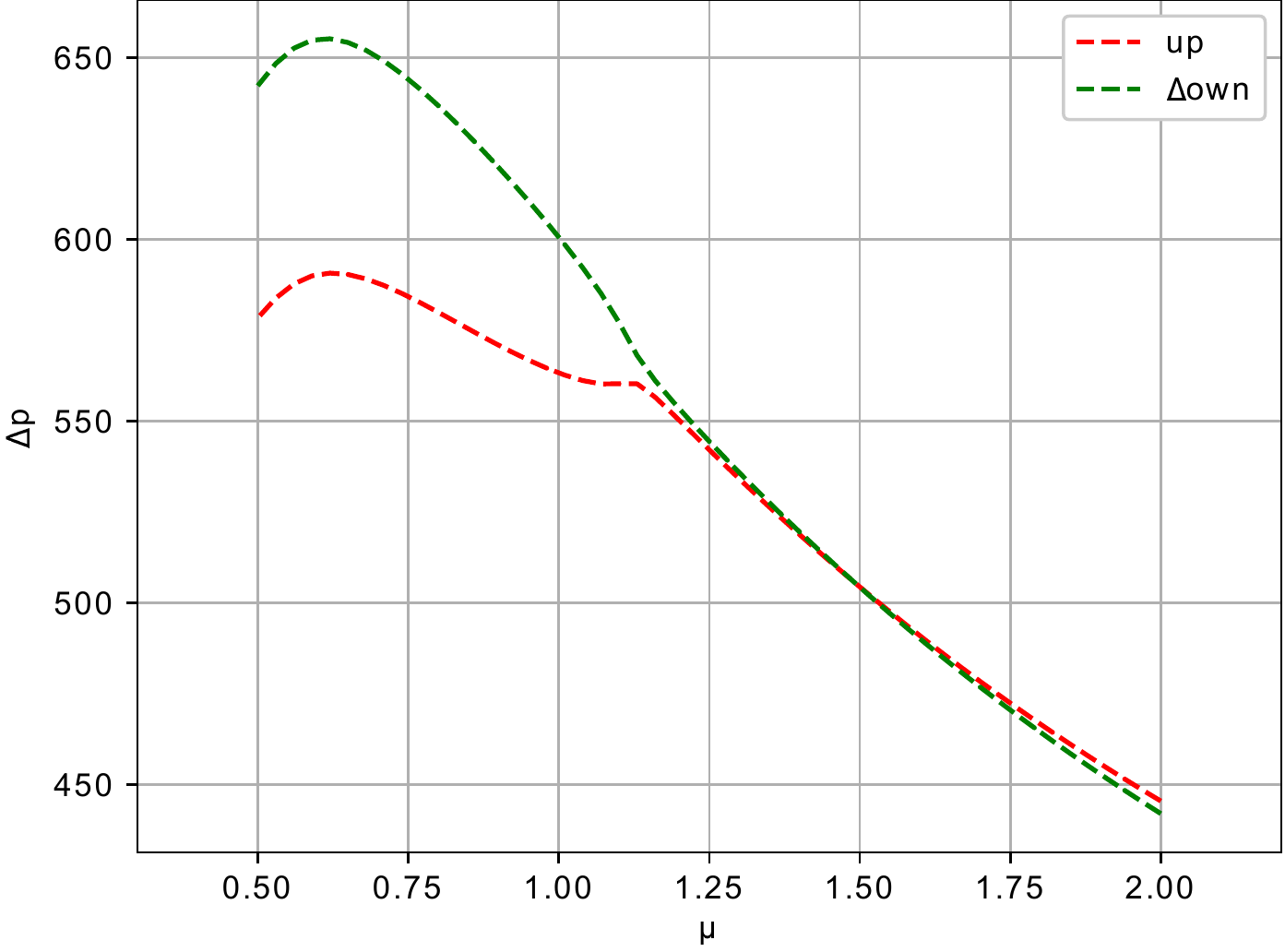} }}%
    \caption{}%
\end{figure}

The result obtained is justifiable in the light of the analysis of the results of Figure \ref {fig:snapshots-fsi}. In fact, an asymmetric flow downstream of the expansion causes the stress field to be asymmetric as well. In particular, an attempt can be made to trace back the deformations' profile to the gradient of the stresses on the two interfaces of each leaflet. For example, let us consider the pressure drop downstream the expansion, that can be roughly quantified through the following quantities:
$$
\Delta p_{up}=p(\boldsymbol{x}_{1})-p(\boldsymbol{x}_{2}), \quad \text{with} \quad \left\{\begin{array}{l}\boldsymbol{x}_{1}=(4.5,6), \\\boldsymbol{x}_{2}=(6.5,6);\end{array}\right.
$$
$$
\Delta p_{down}=p(\boldsymbol{x}_{3})-p(\boldsymbol{x}_{4}), \quad \text{with} \quad \left\{\begin{array}{l}\boldsymbol{x}_{3}=(4.5,1), \\\boldsymbol{x}_{4}=(6.5,1).\end{array}\right.\\
$$
%
The trend of these indicators as a function of the kinematic viscosity $ \mu $ is shown in Figure \ref {fig:deltaptrend}. Hence, the pressure gradient presents a profile capable of justifying an asymmetry in the field of displacements of the solid phase. A more precise analysis would require calculating the result of the complete stress tensor (which also includes the viscous component linked to the velocity gradient) and integrating on each leaflet's surface.

\begin{figure}[h]
    \subfloat[\label{fig:deltad}Bifurcation diagram for the solid phase.]{{\includegraphics[width=0.5\linewidth]{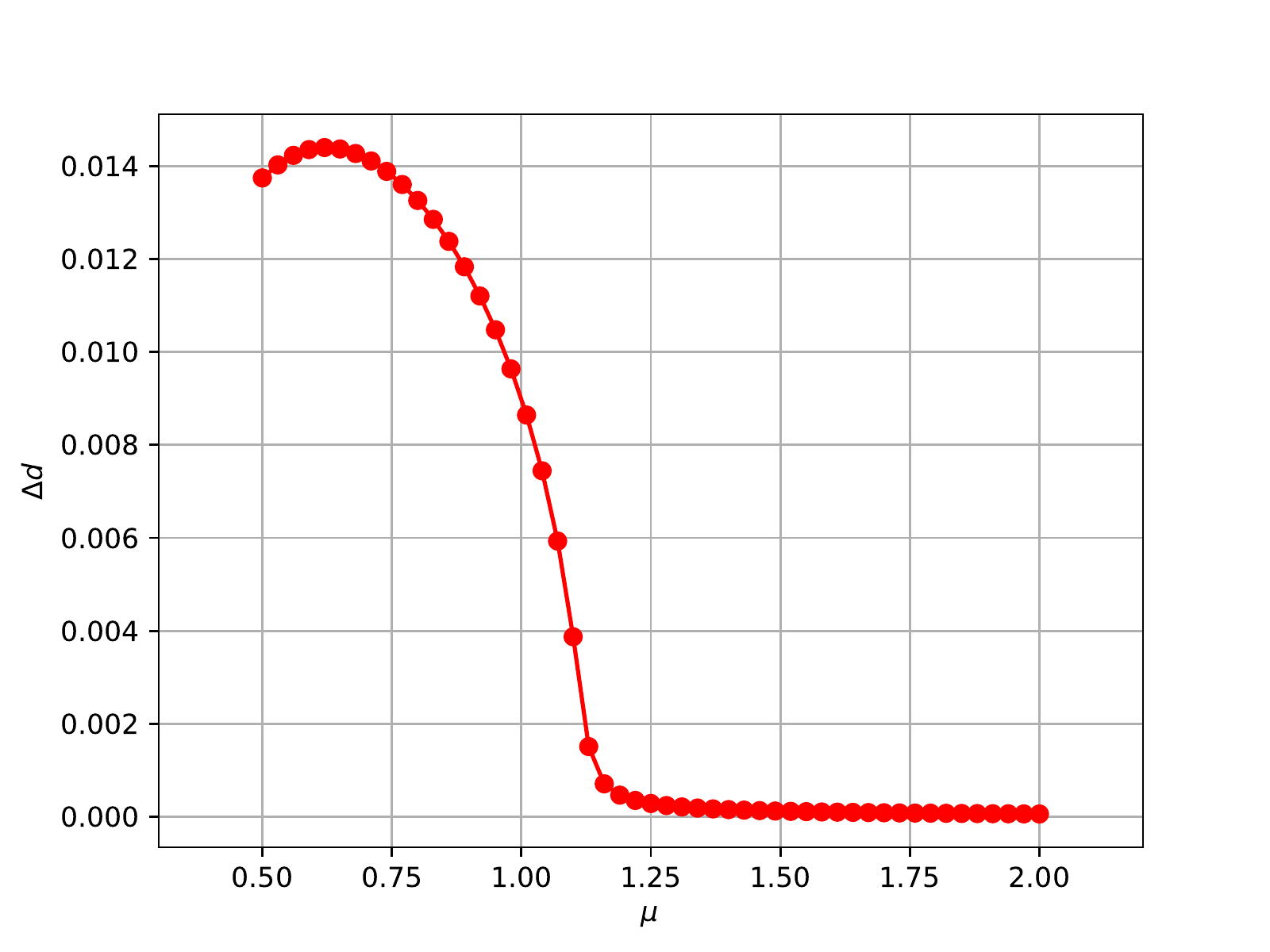}}}
    \subfloat[\label{fig:vmax}Maximum magnitude of the velocity in correspondence to the expansion.]{{\includegraphics[width=0.50\linewidth]{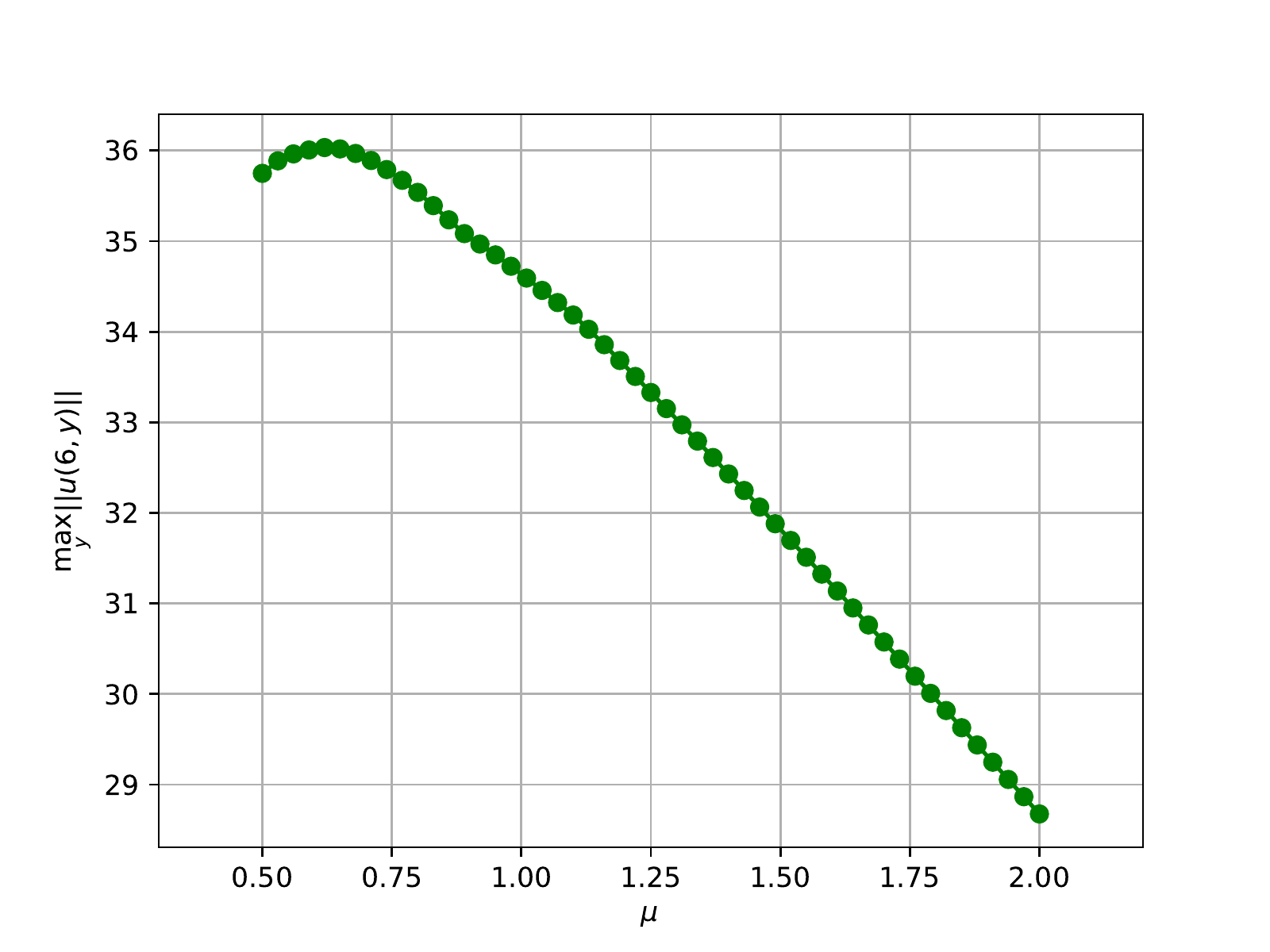}}}
    \caption{}
\end{figure}

In order to retrieve a bifurcation diagram for the solid phase, one can exploit the asymmetric trend of the upper and lower deformations below the critical value $ \mu ^ {*} $. If we indicate with $ \hat {\Omega} _ {1} $ the upper leaflet and with $ \hat {\Omega} _ {2} $ the lower one (Figure \ref{fig:deltap}), we can then measure the difference between the deformations with the indicator $\Delta d= | \max_{\hat{x}\in \hat{\Omega}_{1}} \hat{\boldsymbol{d}}_{s} - \max_{\hat{x}\in \hat{\Omega}_{2}} \hat{\boldsymbol{d}}_{s}|.$
If we plot the trend of this quantity as a function of $ \mu $ we obtain the bifurcation diagram shown in Figure \ref{fig:deltad}.
Therefore, in this case, we arrive at a result that is similar to the bifurcation diagram for the fluid phase (Figure \ref {fig:bifurcation-fsi}).

\begin{figure}[h]
    \centering
    \includegraphics[width=1\linewidth]{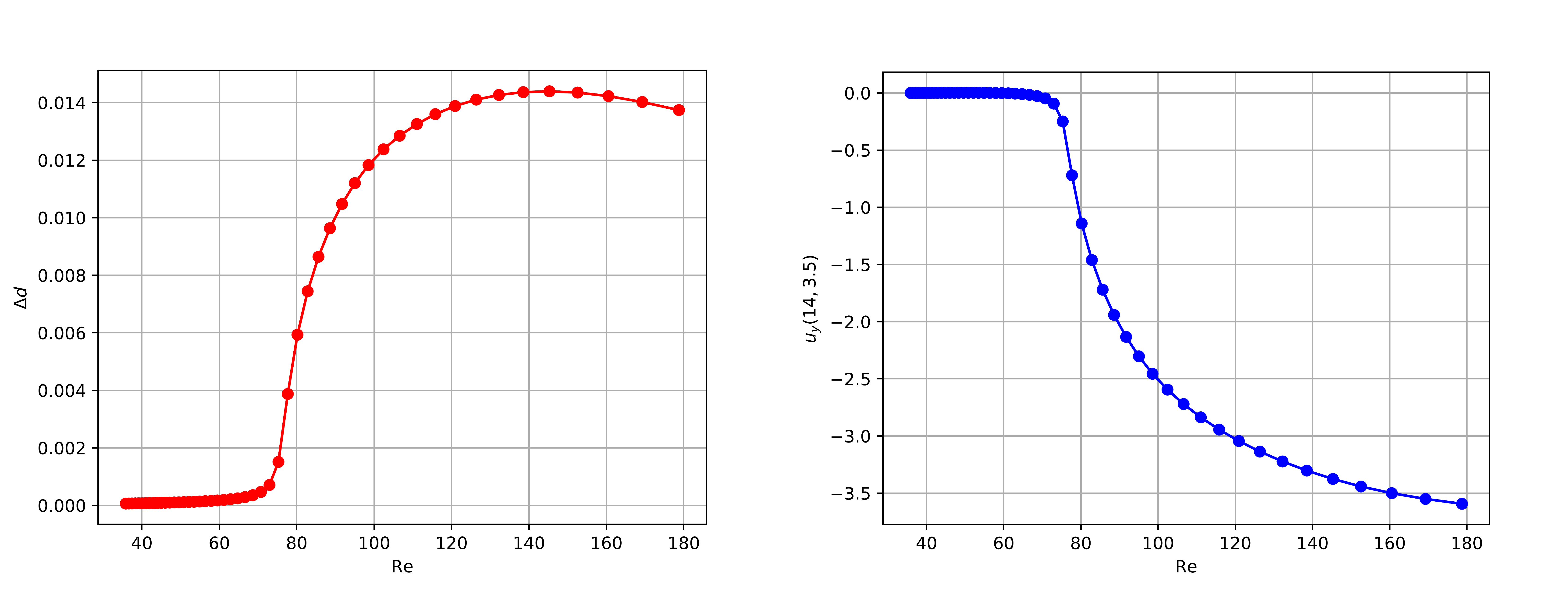}
    \caption{Bifurcation diagrams w.r.t.\ the solid (left) and to the fluid (right) solution's components.}%
    \label{fig:both-diagrams}
\end{figure}

We want to remark that developing a bifurcation analysis as a function of the Reynolds number is possible. We can associate a Reynolds number to each viscosity value if we consider the expansion orifice's width as the characteristic length. This can be done using the maximum characteristic velocity at the inlet of the channel's width section. As anticipated at the end of Section \ref {sec:rigid_structure}, this velocity is not constant as $ \mu $ varies, and  its trend is shown in Figure \ref{fig:vmax}.
At this point, for each viscosity value, we compute the Reynolds number as $Re=\frac{U L}{\mu}$; and through it, we express the bifurcation diagrams as in Figure \ref{fig:both-diagrams}, where we have $Re^{*}\approx 75$.
\subsubsection{Reduced order model results}%
\label{sub:reduction_model}
Now we will present the results of the reduction phase for the problem above.
To this end, we used the snapshots coming from the continuation method presented in the previous section. Using the notation of Section \ref{podgalerkin}, this means that $\mathbb{P}_{train}$ is the discretized interval [0.5,2], and $N_{train}=51$.
We then processed the snapshots to extract a reduced basis for each of the fields involved in the FSI problem. In this regard, we remind that there are 6 unknown fields\footnote{In particular they are $\{\mathbf{u}_{f},p,\boldsymbol{d}_{s},\boldsymbol{d}_{f},\boldsymbol{\lambda}_{u},\boldsymbol{\lambda}_{d}$\}, however only the first three have a physical meaning. The others are the Lagrange multipliers introduced to impose boundary conditions and the mesh displacement which is needed for the ALE formulation.} and consequently, for each of them, it was necessary to perform a POD.
However, the reduced velocity space was also expanded using a supremizer enrichment to avoid stability issues \cite{bal}.

\begin{figure}[h]
    \centering
    \includegraphics[width=0.6\linewidth]{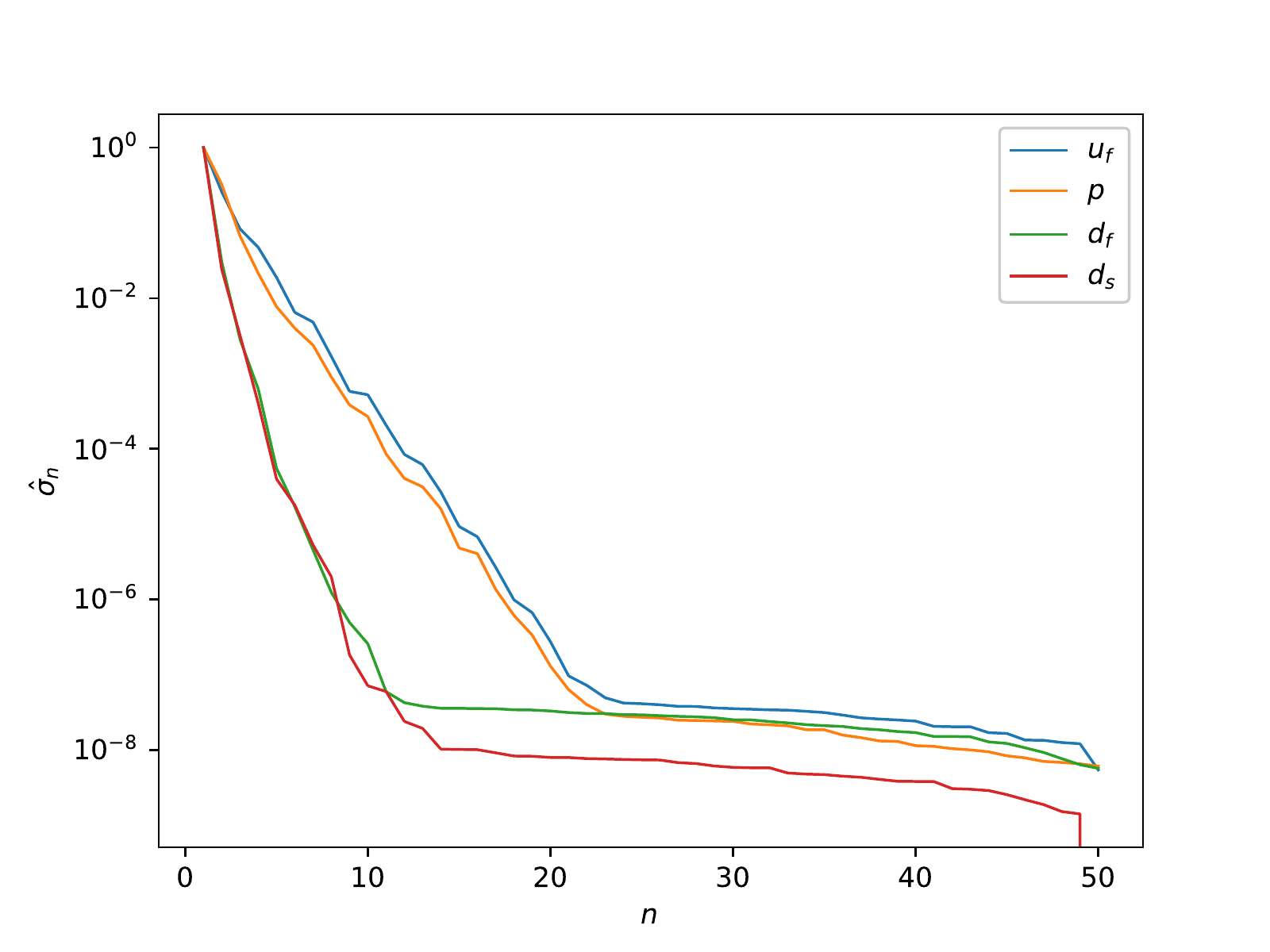}
    \caption{Decay of the normalized singular values for each physical unknown.}%
    \label{fig:energies}
\end{figure}
We show in Figure \ref{fig:energies} the exponential decay of the normalized singular values $\hat{\sigma}_n$ w.r.t.\ the basis cardinality.
In particular, we notice that the two displacement fields decay faster than the pressure and the velocity of the fluid.
In fact, we recall that the problems associated with $ \boldsymbol {d} _ {s} $ and  $ \boldsymbol {d} _ {f}$ are both linear, and this justifies the lower complexity of their solution manifold.

\begin{figure}[h]
\subfloat[\label{fig:avg-projection}Average projection error in logarithmic scale.]{{\includegraphics[width=0.5\linewidth]{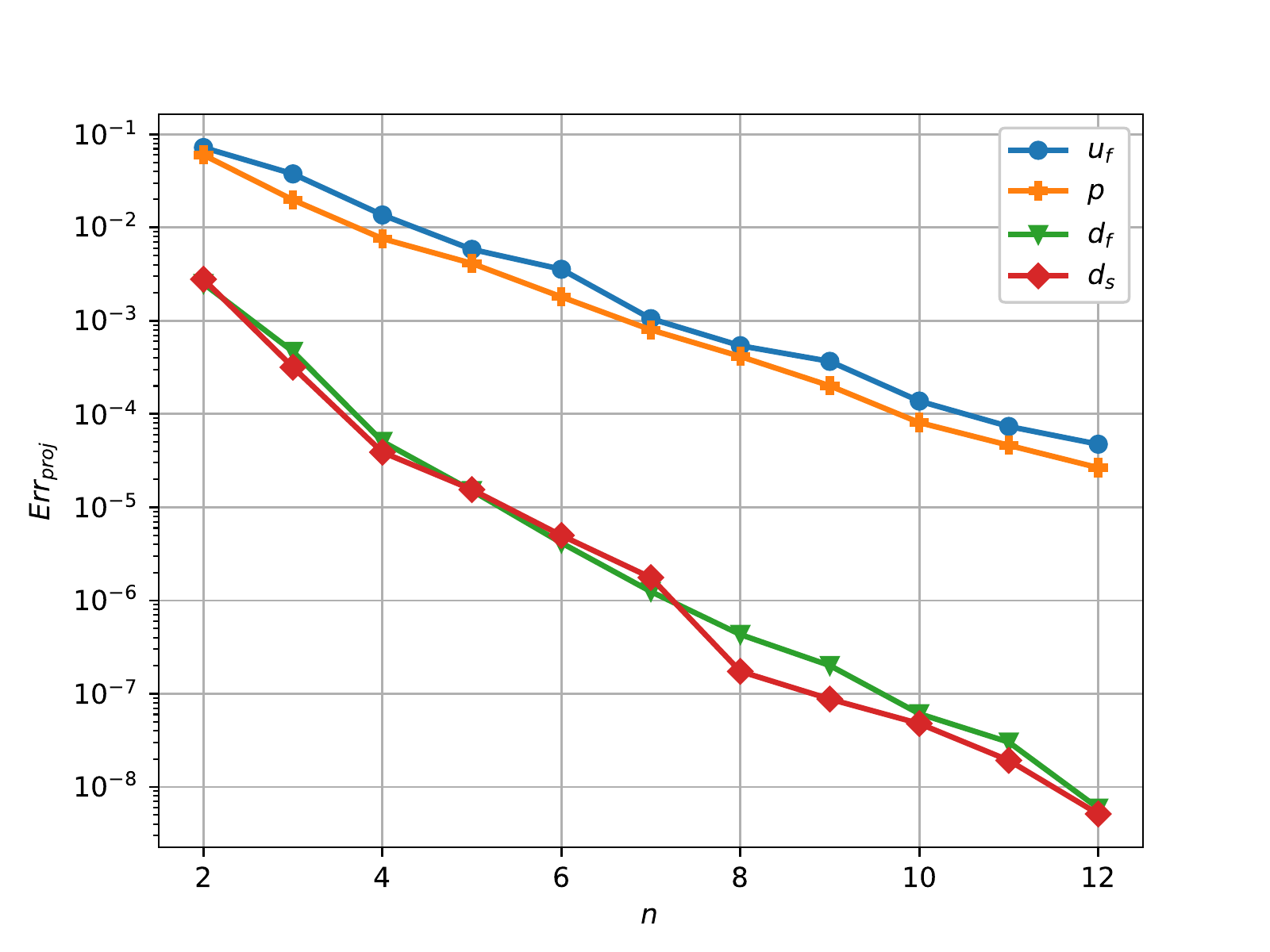}}}
    \subfloat[\label{fig:avg-reduction}Average reconstruction error in logarithmic scale.]{{\includegraphics[width=0.5\linewidth]{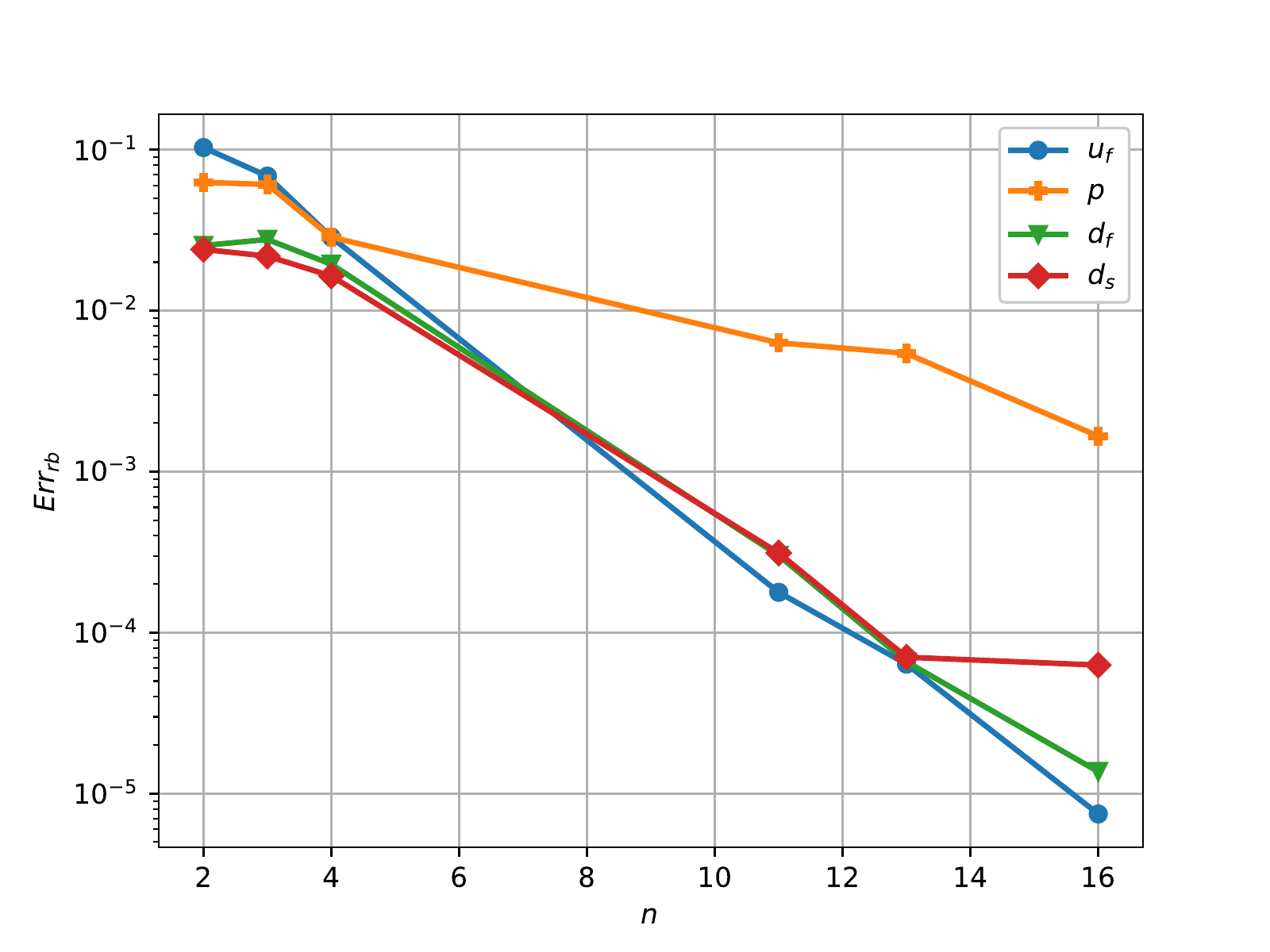}}}
\caption{}
\end{figure}

A valuable tool for evaluating the effectiveness of the basis extracted from the POD is the projection errors' trend as a function of the basis cardinality.
Since we are interested in reconstructing the solution for all the values of $ \mu \in \mathbb {P} _ {train}$, we consider the following average error:
\begin{equation}
    Err_{proj}=\frac{1}{N_{train}} \sum_{j=1}^{N_{train}} \frac{\| X_{h}(\mu_{j})-P_{n}(X_{h}(\mu_{j}))\|_{\mathbb{X}}}{\|X_{h}(\mu_{j})\|_{\mathbb{X}}},
\end{equation}
with $P_{n}: \mathbb{X}_{h} \rightarrow \mathbb{X}_{rb}$ being the projection operator of the high-fidelity solution onto the reduced basis space.
The trend of this indicator for the different fields is shown in Figure \ref{fig:avg-projection}. It should be noted that the projection error represents the \textit {best approximation error}, so one cannot obtain an error lower than that (Cèa's Lemma).

We have analyzed the average error trend for some basis cardinalities in the case of solutions obtained through the Reduced Basis method. In this case, the difference between the high fidelity solutions and the reduced ones was measured as follows:

\begin{equation}
    Err_{rb}=\frac{1}{N_{train}} \sum_{j=1}^{N_{train}} \frac{\| X_{h}(\mu_{j})-X_{rb}(\mu_{j})\|_{\mathbb{X}}}{\|X_{h}(\mu_{j})\|_{\mathbb{X}}}.
\end{equation}

In Figure \ref{fig:avg-reduction} we observe that, as expected, the errors are higher than their direct projection counterparts.
In fact, the nonlinearity of the system, which makes its resolution dependent on the Newton method and the chosen initial guess, and its singularity at the bifurcation point, ensure that the error decay does not have a monotone behavior; however, this decay remains exponential. Furthermore, the pressure field is the one with the highest reconstruction errors, while the velocity field can be reconstructed with the same accuracy as the  displacements.

\begin{figure}[h]
    \centering
    \includegraphics[width=0.7\linewidth]{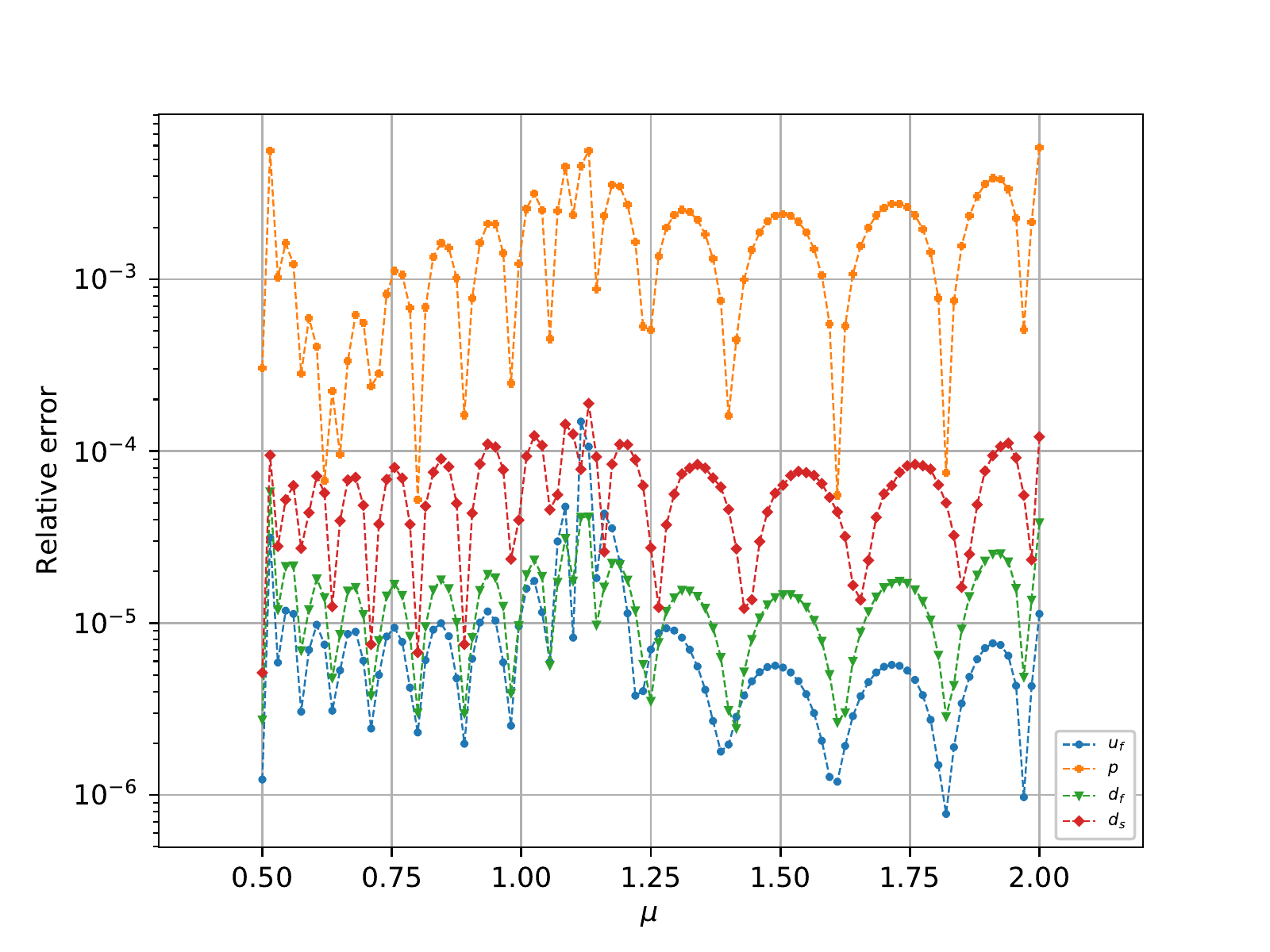}
    \caption{Relative reconstruction error in logarithmic scale with $N_{rb}=16$ for all of the physical variables over the parameter interval [0.5,2].}%
    \label{fig:error}
\end{figure}
At this point, we have decided to work with $ N_ {rb} = $ 16 reduced basis for each one of the fields.
In this regard, to correctly test the effectiveness of reconstruction by the reduced basis, we have refined the discretization of the $ \mathbb {P}_{train} $ parameter interval by considering $ 101 $ equally spaced points.
Following the implementation of the reduced continuation method (see Algorithm \ref {algorithm-reduced}), the reduced solutions for the lower stable branch are obtained.  We then compared them with the full order counterpart through the relative errors for each field taken into consideration (see Figure \ref{fig:error}).

We observe that these errors reflect the behavior already highlighted through the average errors (Figure \ref{fig:avg-reduction}), therefore the pressure is reconstructed less accurately than the other fields.
Moreover, we can observe a slight peak for the error trend in correspondence with the bifurcation value and the low viscosity regime.
Thus it seems that the presence of the structure somehow helps the regularity of the parametric manifold, allowing for a better RB approximation of the bifurcating phenomena w.r.t. the one for the rigid case \cite{Nonino}.

\begin{figure}[h]
\centering
\subfloat[$ \delta \mathbf{u}_{rel}$]{%
  \includegraphics[width=1\textwidth]{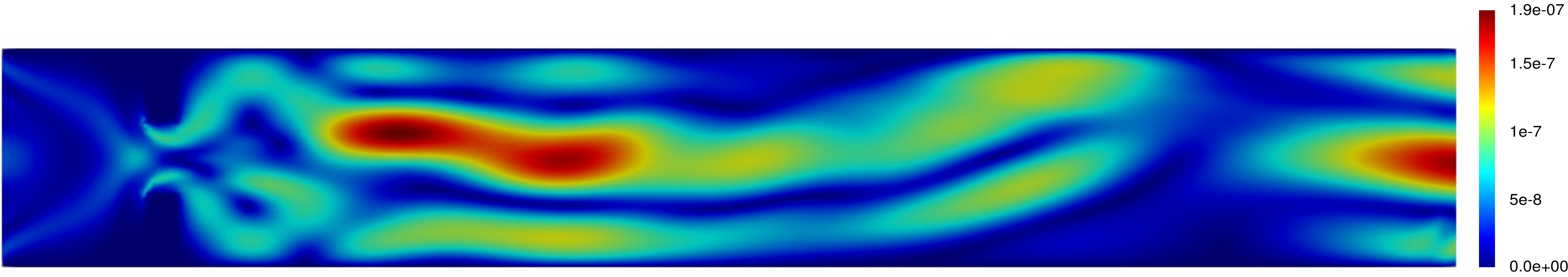}%
  }\par
\subfloat[$\delta p_{rel}$]{%
  \includegraphics[width=1\textwidth]{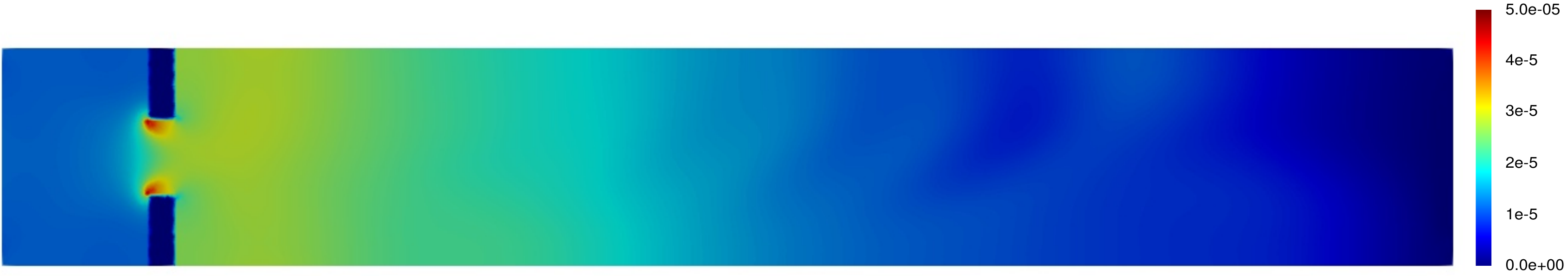}%
  }\par
  \subfloat[$(\delta \boldsymbol{d}_{f,rel},\delta \boldsymbol{d}_{s,rel})$]{%
  \includegraphics[width=1\textwidth]{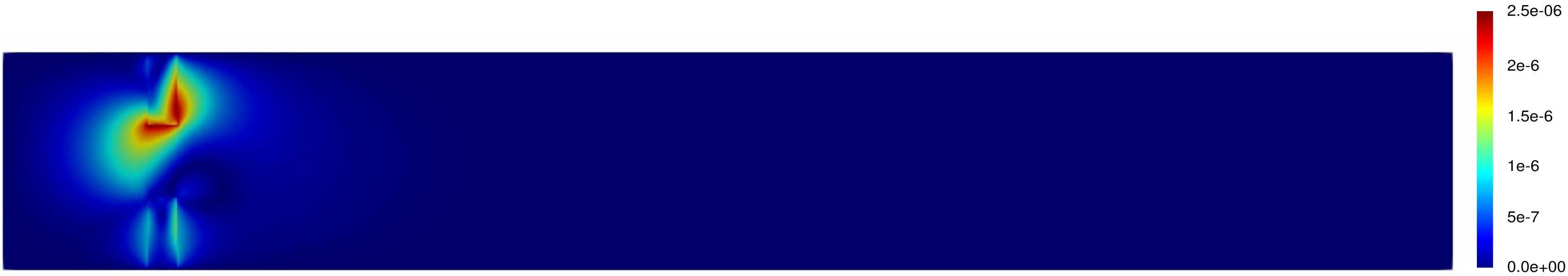}%
  }
\caption{Relative error fields for $\mu=0.5$ : for what concerns the displacement, we recall that that $\delta \boldsymbol{d}_{f,rel}$ (respectively $\delta \boldsymbol{d}_{s,rel}$)  is defined only in $\hat{\Omega}_{f}$ (respectively $\hat{\Omega}_{s}$.}
\label{fig:delta-fields}
\end{figure}
For some application it may be important also to understand the effectiveness of the reconstruction as a function of the spatial coordinate $\hat{\boldsymbol{x}}$. For this reason, we have also computed the component-wise relative error field, which assumes the following expression:
\begin{equation}
    \delta X_{rel}=\frac{|X_{h}(\mu_{j})-X_{rb}(\mu_{j})|}{{\|X_{h}(\mu_{j})\|_{\mathbb{X}}}}.
\end{equation}
%
Figure \ref{fig:delta-fields} shows these relative error associated with the four fields for $\mu=0.5$ . We can observe how each field presents criticalities in different regions of the domain. In particular, the two displacement fields are reconstructed with less efficiency in correspondence with the $ \hat{\Gamma} _{FSI} $ interface.
On the other hand, the pressure has the highest reconstruction error between the leaflets.

Finally, the error for the velocity field traces the magnitude of the field itself, with a morphology similar to the one of the deflected jet, also presenting an additional criticality at the outlet of the channel.


\subsection{FSI problem with the nonlinear Saint Venant-Kirchoff model}%
\label{sec:sectionSVK}

As we have seen in the previous section, introducing a structure with a linear constitutive law does not seem to modify the phenomenon of bifurcation undergone by the fluid phase.

However, it is necessary to consider the fact that the previous results were obtained with a high structural stiffness, which results in a small-scale deformation field. Precisely because of this amplitude, it is reasonable to expect that the bifurcation point is only mildly affected by the parametric interval taken into consideration.

Therefore, now we want to reduce the structural stiffness and analyze the effect this has on the bifurcation. At the same time we will modify the solid phase's constitutive relation, exploiting a nonlinear hyperelastic model. This is done because, for many materials, the linear elastic model fails to describe high strains' behavior accurately. The adjective hyperelastic in this context refers to the fact that the relation between stress and strain can be obtained by deriving a scalar quantity known as the strain energy density function. The simplest hyperelastic model is the \textit{Saint Venant-Kirchoff} (SVK) one, which represents a straightforward extension of the linear elastic model.

If we define the Green-Lagrange strain tensor as
\begin{equation}
    \mathbf{E}=\frac{1}{2}(\mathbf{F}^{T}\mathbf{F}-\mathbf{I})=\frac{1}{2}( \nabla \boldsymbol{d}_{s}+ \nabla ^{T} \boldsymbol{d}_{s}+ \nabla^{T} \boldsymbol{d}_{s} \nabla \boldsymbol{d}_{s}),
\end{equation}

then the strain energy turns out to be a function of the deformation gradient $\mathbf{F}$ through an explicit dependence on $\mathbf{E}$, that is:
\begin{equation}
    \Psi(\mathbf{F})=\mu_{s} \mathbf{E}:\mathbf{E}+ \frac{\lambda}{2}(\operatorname{tr}(\mathbf{E}))^{2},
\end{equation}
where $\lambda$ and $\mu_{s}$ are the Lamé's constants.
Finally, to find the constitutive relation that expresses the Piola tensor is given by
\begin{equation}
    \frac{\partial \Psi}{\partial \mathbf{F}} =\mathbf{P}(\mathbf{F})= \mathbf{F}(2 \mu_{s}\mathbf{E} + \lambda\operatorname{tr}(\mathbf{E}) \mathbf{I}).
\end{equation}
One can replace the two Lamé's constants by another pair of independent parameters, whose values characterize the material under consideration. A very popular choice consists in using the Young's modulus ($E$) and the Poisson's ratio ($\nu$), such that we can write
\begin{equation}
    \mu_{s}= \frac{E}{2(1 + \nu)},\text{ and }  \lambda= \frac{E \nu}{(1+\nu)(1-2 \nu)}.
\end{equation}
In the following we will keep the Poisson's ratio fixed to $\nu=0.21$ , while we will consider different values for the Young’s modulus $E$ (see Table \ref{table6}).

\begin{figure}[h]
    \centering
    \includegraphics[width=0.7\linewidth]{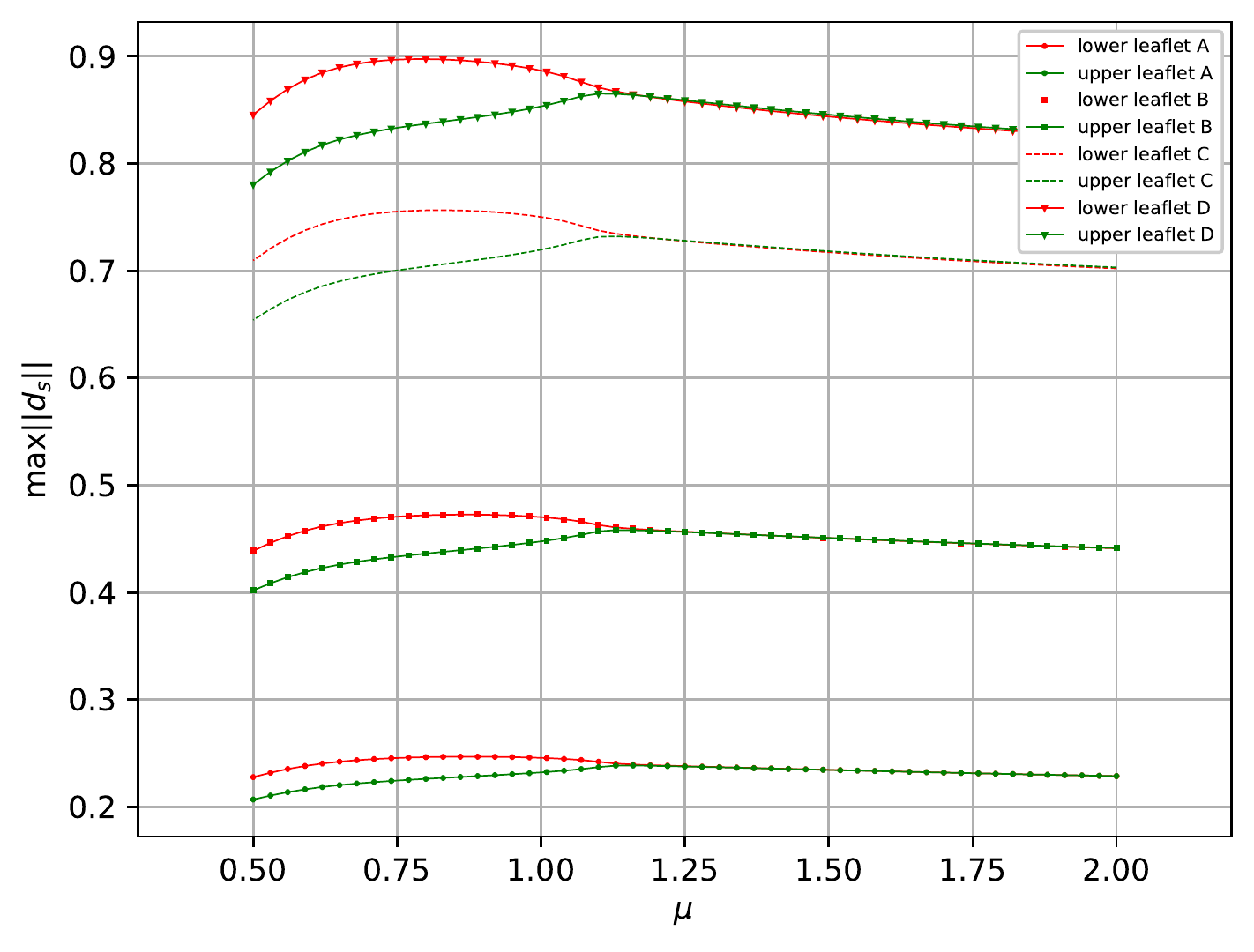}
    \caption{Maximum $\boldsymbol{d}_{s}$ for the different values of $E$ under consideration.}%
    \label{fig:comparison-nonlin-disp}
\end{figure}

\begin {table} [H]
\centering
\begin {tabular} {lc}
    \toprule
    \textbf {Case} & E  \\
    \midrule
    A) & 1.87e5\\
    B) & 9.34e4\\
    C) & 5.60e4\\
    D) & 4.67e4\\
    \bottomrule
\end {tabular}
\caption{Values of the Young’s modulus used to model different stiffness conditions.}
\label{table6}
\end{table}
\begin{figure}[h]
    \centering
    \includegraphics[width=1\linewidth]{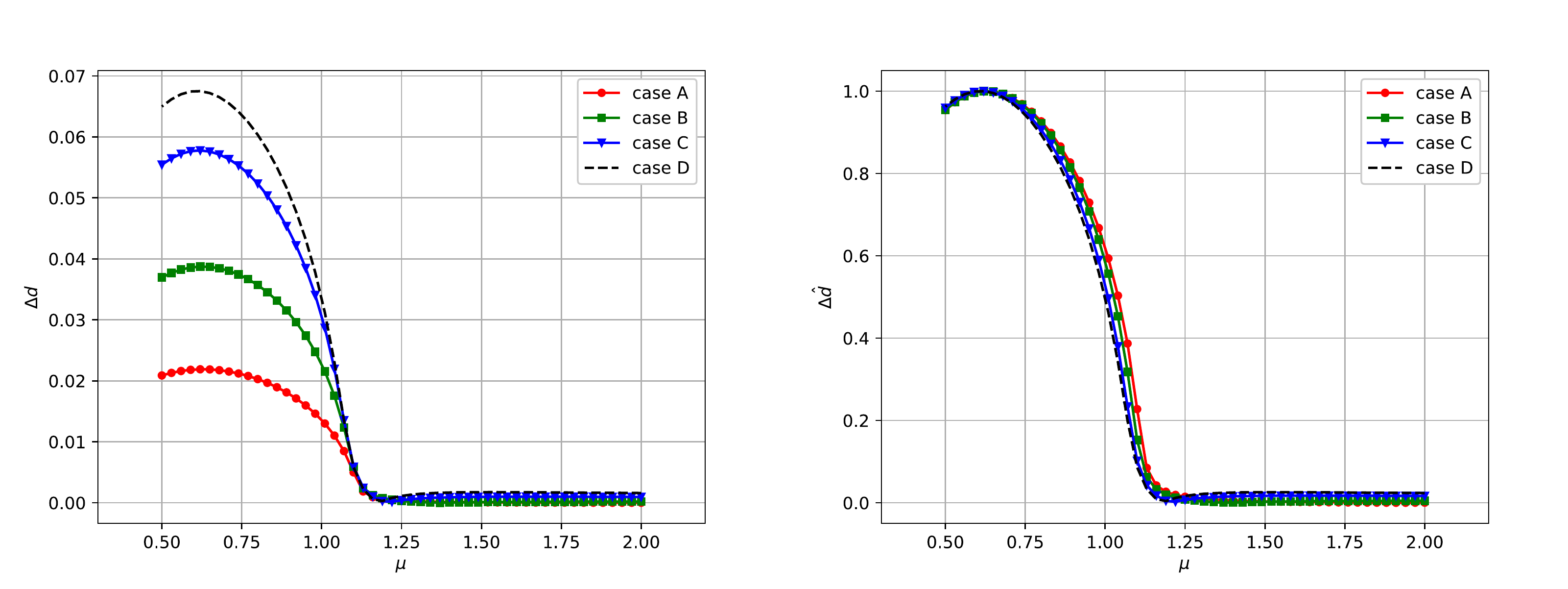}
    \caption{Bifurcation diagrams for the solid phase obtained through $\Delta d$ (left) and $\Delta \hat{d}$ (right). }%
    \label{fig:deltad-comparison}
\end{figure}
\begin{figure}[h]
    \centering
    \includegraphics[width=1\linewidth]{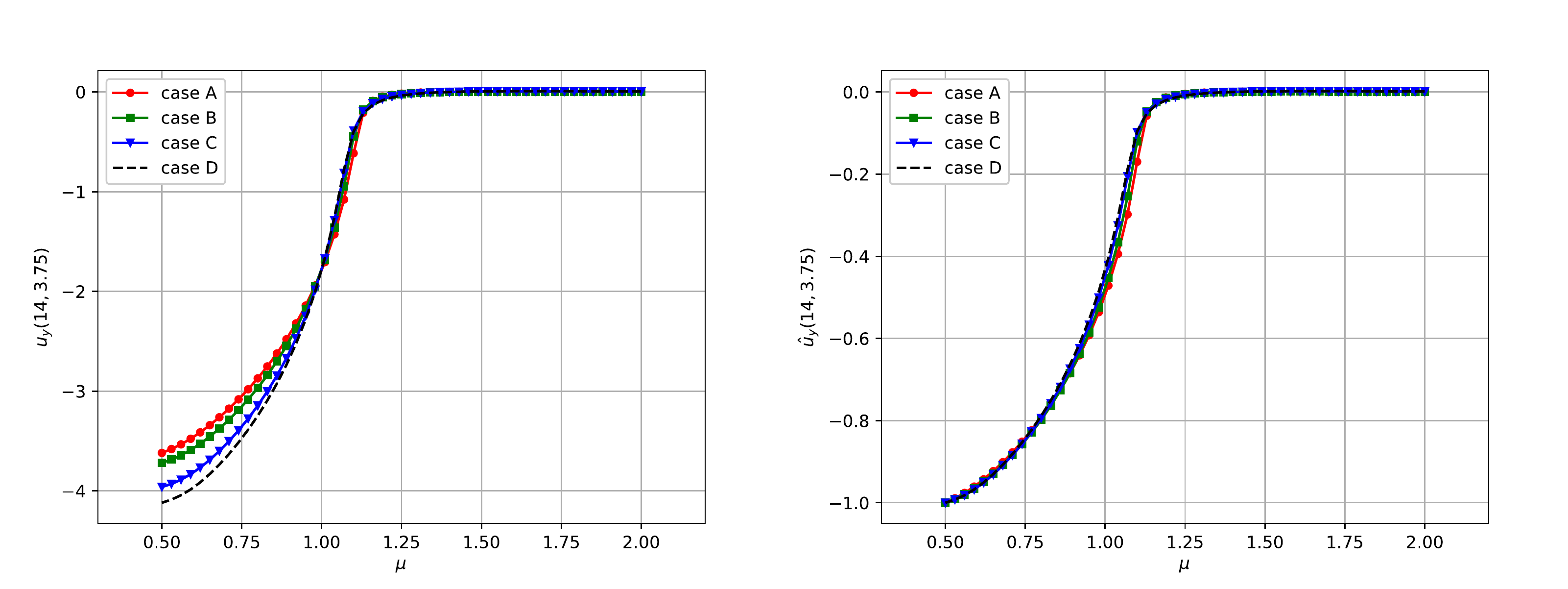}
    \caption{Bifurcation diagrams for the fluid phase obtained through $ u_ {y} (\underline {\boldsymbol {x}})$ (left) and $ \hat{u}_ {y} (\underline {\boldsymbol {x}}) $ (right). }%
    \label{fig:deltau-comparison}
\end{figure}

First of all, we can compare the maximum displacements of the upper and lower leaflets for the different cases, as shown in Figure \ref {fig:comparison-nonlin-disp}.

As can be easily observed, we have reduced the structural stiffness until reaching deformations of the same order as the domain's dimensions.
However, the analysis carried out in Figure \ref{fig:comparison-nonlin-disp} does not yet allow us to understand if the bifurcation behavior is only rescaled due to the differences in the values of the Young's modulus or if there are differences in response depending on the case considered.
Because of the former consideration, it makes sense to use as an indicator the difference between the maximum deformations of the upper and lower leaflets $ \Delta d $, which was used in the previous section to compute a bifurcation diagram for the solid phase.
The result is shown in Figure \ref{fig:deltad-comparison}. In the same figure, we also present the trend of its dimensionless counterpart $ \Delta \hat {d} $, defined as
\begin{equation}
    \Delta \hat{d} \doteq \frac{ \Delta d }{ \max_{{\mu \in \mathbb{P}_{h}}} \Delta d}.
\end{equation}
We can observe, it can be observed, through the trend of the normalized indicator $ \Delta \hat{d} $, that the bifurcation behavior of the structure is only scaled by the stiffness (Young's modulus).
A final analysis with this respect consists in comparing the various bifurcation diagrams obtained for the fluid phase through the scalar output used previously, namely $ u_ {y} (\underline {\boldsymbol {x}}) $, with $ \underline {\boldsymbol {x}} = ( 14,3.75)$. Also in this case, the normalization with respect to the maximum value allows us to conclude that even for the fluid phase the behavior is simply scaled with respect to the value of $ E $ (see Figure \ref{fig:deltau-comparison}).
\subsection{Model comparison for the FSI problem}%
\label{sec:comparison}
Finally, we compare the results for the different models considered. We recall that the differences concern the modelling of the leaflets for which we used a rigid body model (Section \ref{sec:rigid_structure}), a linear elastic model (Section \ref{sec:fsi_problem_with_linear_elastisticy}), and the Saint Venant-Kirchoff nonlinear model (Section \ref{sec:sectionSVK}). However, a rigorous comparison between the linear and the nonlinear model requires that the same values of the parameters in the constitutive relations are used. Therefore, we decided to add a test case (A) for the linear model, for which $ E = 4.67e4 $ and $ \nu = 0.21 $ (corresponding to the test case (D) of the former section). In particular, we recall that these values result in a structural deformation of the same order of magnitude of the domain's dimensions.
Table \ref{table-comparison} shows the values of the parameters for the different models compared in this section.

\begin{figure}[h]
\subfloat[\label{fig:comparison-lin-nonlin}Maximum displacement of the leaflets for both the SVK and the linear models.]{{\includegraphics[width=0.48\linewidth]{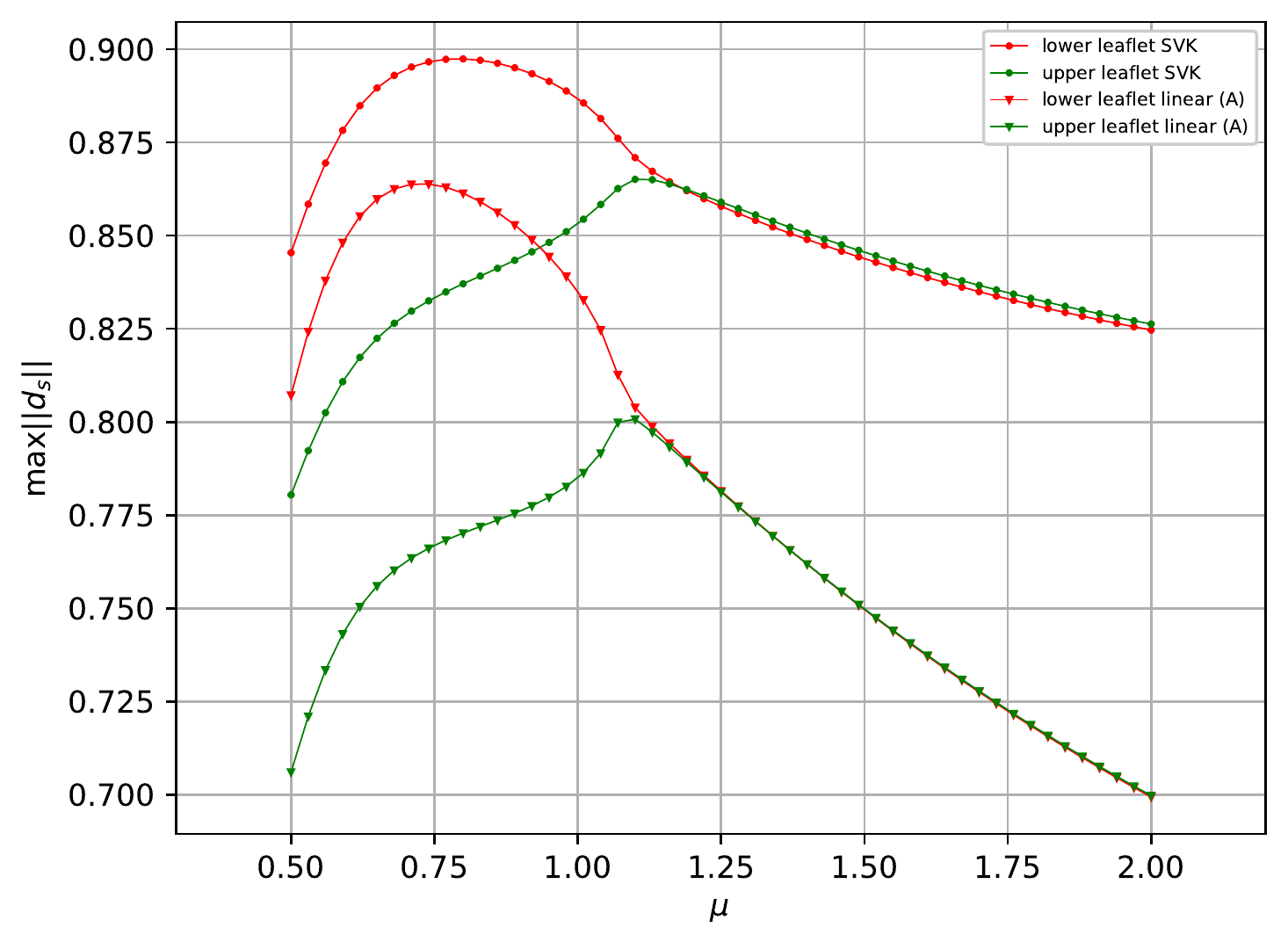}}}
\subfloat[\label{fig:deltad-lin-nonlin} Difference between the displacement of the upper and lower leaflets, as measured through $\Delta d$, for both the SVK and the linear models.]{{\includegraphics[width=0.52\linewidth]{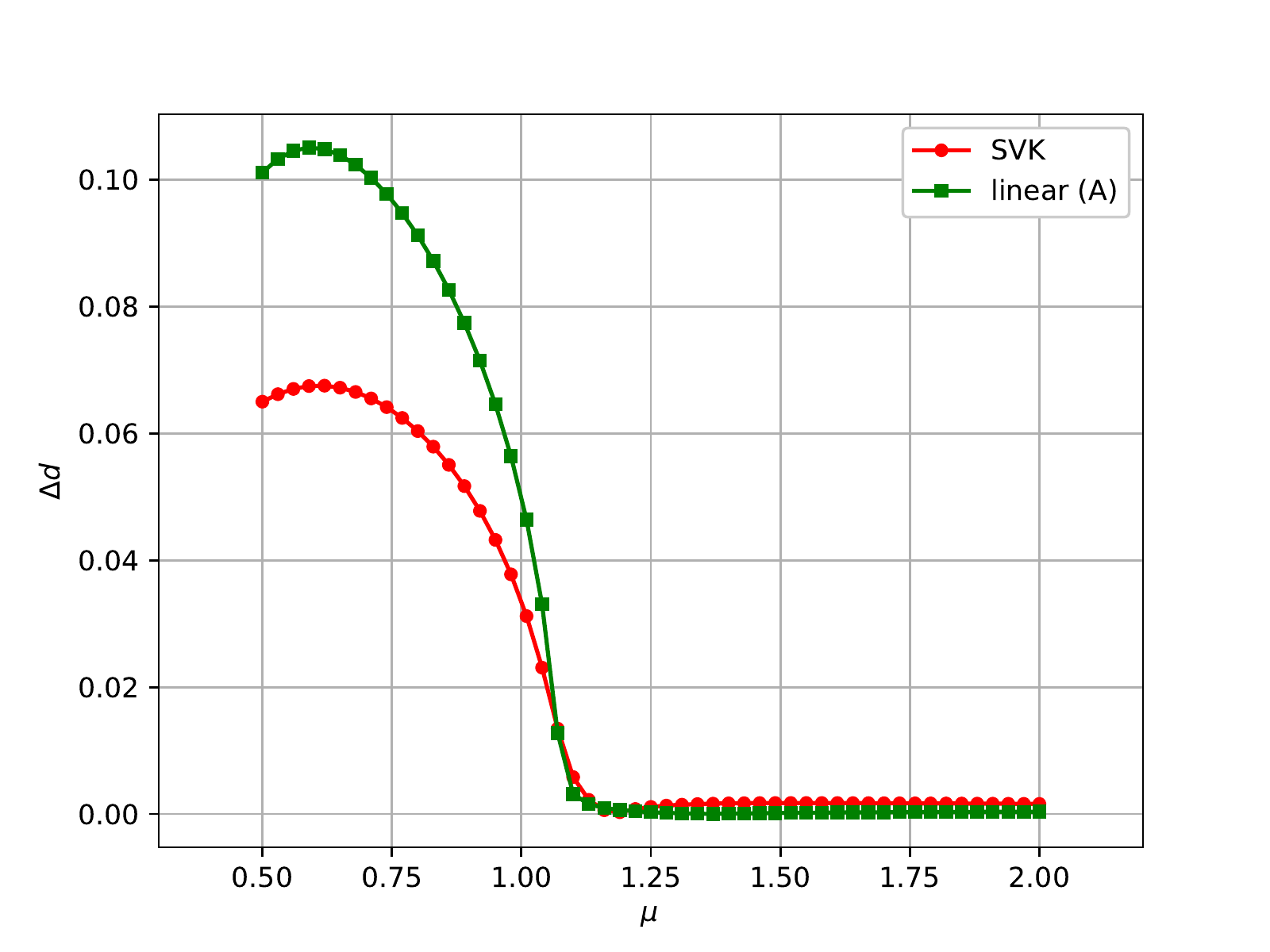}}}
\end{figure}

\begin{figure}[H]
    \centering
    \includegraphics[width=0.6\linewidth]{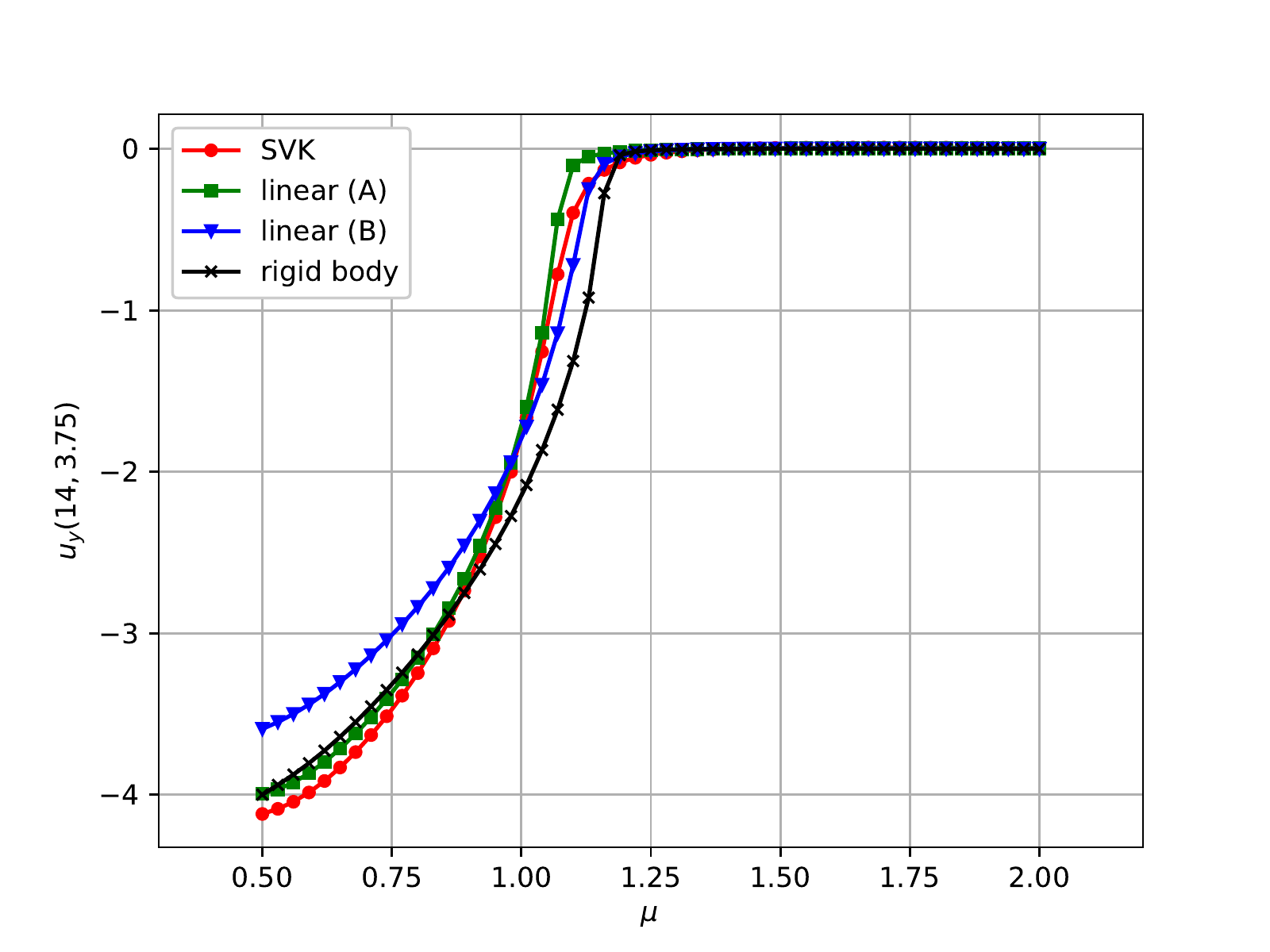}
    \caption{Comparison of fluid phase bifurcation diagram for the different models.}%
    \label{fig:vy-comparison-lin-nonlin}
\end{figure}

Let us start the comparison between the linear model (A) and the SVK model.
In Figure \ref{fig:comparison-lin-nonlin} we examine the maximum displacements of the two leaflets. The first observation concerns the fact that the SVK model corresponds to more significant deformations.
\begin {table} [h]
\centering
\begin {tabular} {lcc}
    \toprule
    \textbf {Case} & $E$  & $\nu$   \\
    \midrule
    linear (A) & 4.67e4&0.21\\
    linear (B) & 2.89e5&0.44\\
    SVK & 4.67e4&0.21\\
    \bottomrule
\end {tabular}
\caption{Values of the Young's modulus and the Poisson's ratio used for the test cases.}
\label{table-comparison}
\end{table}
Furthermore, we can observe an asymmetry in the deformation of the two leaflets for the SVK model even before of the bifurcation point, while this does not occur for the linear counterpart.
Another interesting fact is the behavior in the post-bifurcation regime. In fact, the linear model predicts that the asymmetry of the two leaflets' behavior is more accentuated than the counterpart represented by the SVK model. This observation can be quantified by comparing the diagrams related to the scalar output $\Delta d$ (see Figure \ref{fig:deltad-lin-nonlin}).
Therefore, the linear model, to which a field of minor deformation is associated, responds to the bifurcation in a more accentuated way, due to its difficulty to model large displacements.
Finally, to investigate how this difference affects the behavior of the fluid phase, we compare the bifurcation diagrams related to the vertical velocity $u_{y}(14,3.75)$.
The result is presented in Figure \ref{fig:vy-comparison-lin-nonlin}.
As expected, changing the structure's constitutive relation has important effects on the fluid flow field, and this also affects the bifurcation point. In particular, the linear model (B) (which corresponds to a more stiff structure) seems to share the bifurcation point with the rigid structure.
However, in the post-bifurcation regime this model gives rise to a vertical flow of lower intensity compared to the other cases examined.
The decrease in structural stiffness in the linear case (A) leads to a delay of the bifurcation with respect to the kinematic viscosity. However, after the bifurcation, this model results in a vertical velocity that closely approximates the rigid structure's case.

Furthermore, the SVK model delays the bifurcation with respect to the case of the rigid structure but anticipates it with respect to the linear (A) counterpart.

Finally, we can note that the SVK model differs more significantly from the linear model precisely in correspondence with the bifurcation, while in the post-bifurcation regime, there is a similar behavior.

\section{Conclusions and perspectives}
\addcontentsline{toc}{chapter}{Conclusions}
The goal of this work was to analyze a complex phenomenon of bifurcation in the multiphysics context of fluid-structure interaction problems. In particular, we have adopted a very general algorithm for the study of the phenomenon at hand. Although the Coand\u{a} effect and the consequent investigation of the bifurcation it causes have already been addressed in previous works, we are not aware of studies on the effect that the coupling of the fluid with an elastic solid can have on the bifurcation.

In line with our main objective, we compared the rigid structure case with that of elastic structure. In a first test case, conducted using a linear constitutive relationship for the solid phase, we observed that the structure passively undergoes a bifurcation, as shown by the leaflets’ displacements. Furthermore, the presence of the structure seems to rescale the bifurcation diagram. On this test case, we also conducted a reconstruction using the Reduced Basis method. The peculiarity of this procedure concerns that the setting of the problem faced is very complex: we are dealing with a nonlinear system with seven unknown fields, and which undergo a bifurcation in the parametric range of interest.
Nevertheless, we were able to extract a basis able to accurately approximate the solution manifold, as demonstrated by the trend of construction errors.

We then decided to deepen our investigation by considering large scale deformations. For consistency with the physical reality of the phenomenon, this was done by introducing a nonlinear constitutive relationship for the structure, i.e. the Saint Venant-Kirchoff model.

In this case, we have observed that modifying the structural stiffness does not affect the bifurcation point's position. However, the comparison between the case of the rigid structure and that of the linear model showed that the solid model actually affects the bifurcation point. We have, in fact, noticed that the introduction of a solid with a nonlinear constitutive law seems to delay the bifurcation, and this delay is amplified in the linear case with large scale deformations.

During the conclusion of this work, we have identified many interesting topics to be deepened.
In the present work, we have not investigated the stability of the solutions. In particular, we intend to develop an eigenvalue analysis that makes it possible to have a quantitative indicator for the bifurcation points; this will allow us to deepen the comparison between the presented test cases.

Moreover,  we have used a monolithic approach for both the high fidelity problem and the reduced order one. A possible extension concerns the use of partitioned schemes based on the Chorin–Temam decomposition. This approach could also be integrated with a semi-implicit scheme for the treatment of the coupling conditions between fluid and solid \cite{Nonino}. Furthermore, the reduced method we used depends on the number of degrees of freedom of the high fidelity problem. Consequently, one possible development concerns the implementation of hyper-reduction strategies that guarantee a more efficient decoupling between the offline and online phases, in order to obtain a more significant speedup.

An interesting perspective results in the integration of our framework with machine learning techniques. In particular, \cite{Pichi}  presents an approach for the efficient detection of the bifurcation points, based on a neural network that exploits the Proper orthogonal decomposition (POD-NN). This strategy can be used to obtain an efficient decoupling between offline and online stages, as an alternative to classical hyper-reduction strategies. However, there are also other approaches of possible investigation, such as the physics informed neural networks (PINNs).

Finally, it is interesting to modify the physics of the problem we have focused on. This can be done for example by considering a setting where the bifurcation is governed by the nonlinearity of the constitutive relationship for the solid problem. For example, we could exploit the buckling results shown for certain nonlinear hyperelastic models. A further generalization results in the investigation of bifurcation phenomena for compressible fluid dynamics. This last problem is considerably complex and of strong interest, as it could also allow to investigate the relationship between bifurcations and turbulence.

\section{Acknoledgments}

We acknowledge the support by European Union Funding for Research and Innovation -- Horizon 2020 Program -- in the framework of European Research Council Executive Agency: Consolidator Grant H2020 ERC CoG 2015 AROMA-CFD project 681447 "Advanced Reduced Order Methods with Applications in Computational Fluid Dynamics". We also acknowledge the PRIN 2017 "Numerical Analysis for Full and Reduced Order Methods for the efficient and accurate solution of complex systems governed by Partial Differential Equations" (NA-FROM-PDEs), and the "GO for IT" program within a CRUI fund for the project "Reduced order method for nonlinear PDEs enhanced by machine learning".

We would also like to show our gratitude to Professor Simona Perotto (Politecnico di Milano) for her useful insights and suggestions which helped shape both the research and the manuscript.


\bibliographystyle{abbrv}
\bibliography{references}

\end{document}